\documentclass[3p,sort & compress]{elsarticle}
\usepackage{amssymb}
\usepackage{latexsym}
\usepackage{amsfonts}
\usepackage{amsthm}
\usepackage{amsbsy}
\usepackage{amsmath}
\usepackage{graphics,graphicx}
\usepackage{subfigure}
\usepackage{multirow,multicol}
\usepackage{bm}
\usepackage{bbm}
\usepackage{color}
\usepackage{float}
\usepackage{calc}
\usepackage{caption}
\usepackage{dsfont}
\usepackage{ifthen}
\usepackage{mathrsfs}
\usepackage[colorlinks,linkcolor=blue,citecolor=blue]{hyperref}
\usepackage{epstopdf}
\usepackage{epsfig}
\usepackage{algorithm}
\usepackage{algorithmic}
\usepackage{pifont}
\usepackage{natbib}
\usepackage{overpic}
\usepackage{psfrag}
\usepackage{rotating}
\usepackage{stmaryrd}
\usepackage{verbatim}
\usepackage{setspace}
\usepackage{tikz}
\usetikzlibrary{fpu}

\newdefinition{remark}{Remark}

\newdefinition{method}{Method}
\newdefinition{example}[theorem]{Example}

\numberwithin{theorem}{section}

\newcommand{\orcid}[1]{\href{https://orcid.org/#1}{\includegraphics[width=8pt]{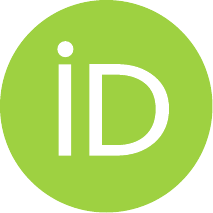}}}
\journal{Journal of \LaTeX\ Templates}
\begin{document}
\begin{frontmatter}
    \title{Fourier heuristic PINNs to solve the biharmonic equations based on its coupled scheme}
    \author[uq]{Yujia Huang} \ead{rainfamilyh@gmail.com}
    \author[Ceyear,SDUInfo]{Xi'an Li\orcid{0000-0002-1509-9328}\corref{cor1}}\ead{lixian9131@163.com}
    \author[acu]{Jinran Wu\orcid{0000-0002-2388-3614}} \ead{ryan.wu@acu.edu.au}
    \cortext[cor1]{Corresponding author.}
    \address[uq]{The University of Queensland, St Lucia 4072, Australia}
    \address[Ceyear]{Ceyear Technologies Co., Ltd, Qingdao 266555, China}
    \address[SDUInfo]{School of Information Science and Engineering, Shandong University, Qingdao 266237, China}
    \address[acu]{Australian Catholic University, Banyo 4014, Australia}
    
\begin{abstract}
Physics-informed neural networks (PINNs) have been widely utilized for solving a range of partial differential equations (PDEs) in various scientific and engineering disciplines. This paper presents a Fourier heuristic-enhanced PINN (termed FCPINN) designed to address a specific class of biharmonic equations with Dirichlet and Navier boundary conditions. The method achieves this by decomposing the high-order equations into two Poisson equations. FCPINN integrates Fourier spectral theory with a reduced-order formulation for high-order PDEs, significantly improving approximation accuracy and reducing computational complexity. This approach is especially beneficial for problems with intricate boundary constraints and high-dimensional inputs. To assess the effectiveness and robustness of the FCPINN algorithm, we conducted several numerical experiments on both linear and nonlinear biharmonic problems across different Euclidean spaces. The results show that FCPINN provides an optimal trade-off between speed and accuracy for high-order PDEs, surpassing the performance of conventional PINN and deep mixed residual method (MIM) approaches, while also maintaining stability and robustness with varying numbers of hidden layer nodes.

    
    \end{abstract}
\end{frontmatter}

\section{Introduction}\label{sec:01}

Partial differential equations (PDEs) are fundamental to a wide array of problems in science and engineering, spanning disciplines such as physics, chemistry, biology, and image processing. Due to their importance, the investigation of PDE solutions has attracted substantial research interest over the years. As exact solutions to PDEs are rarely accessible, numerical approximation methods have become indispensable for their resolution.

In this study, we concentrate on the numerical solution of the biharmonic equation: 
\begin{equation}\label{eq:biharmonic}
\Delta^2 u(\bm{x}) = f(\bm{x}, u, \Delta u)  ~~\text{in}~ \Omega
\end{equation}
with two different cases of boundary conditions:
\noindent\begin{itemize}
	\item Dirichlet boundary conditions
	\begin{equation}\label{eq:Dirichlet}
	\displaystyle ~u(\bm{x}) = g(\bm{x}) ~~\text{and}~~ \frac{\partial u(\bm{x})}{\partial \vec{n}}= h(\bm{x}) ~~\text{on}~ \partial\Omega,
	\end{equation}
	\item Navier boundary conditions
	\begin{equation}\label{eq:Navier}
	\displaystyle ~u(\bm{x}) = g(\bm{x}) ~~\text{and}~~\Delta u(\bm{x})= k(\bm{x}) ~~\text{on}~ \partial\Omega.
	\end{equation}
\end{itemize}
Consider a domain $\Omega$, which is a polygonal or polyhedral region in Euclidean space $R^d$, with a boundary that is piecewise Lipschitz and satisfies the interior cone condition. Let $f(\bm{x}, u) \in L^2(\Omega)$ represent a known function, and the outward normal component on $\partial\Omega$ is given by ${\partial u}/{\partial \vec{n}}$. Here, $\Delta$ denotes the standard Laplacian operator. The expressions for $\Delta u(\bm{x})$ and $\Delta^2 u(\bm{x})$ are as follows:
\begin{equation*}
\Delta u(\bm{x}) = \sum_{i=1}^{d} \frac{\partial^2 u}{\partial x_i^2} \quad \text{and} \quad \Delta^2 u(\bm{x}) = \sum_{i=1}^{d} \sum_{j=1}^{d} \frac{\partial^4 u}{\partial x_i^2 x_j^2},
\end{equation*}
respectively.

The biharmonic equation is a frequently encountered high-order partial differential equation that arises in fields such as physics and mathematics, particularly within elasticity theory and Stokes flow problems. Applications of this equation include scattered data fitting with thin plate splines \cite{wahba1990spline}, fluid mechanics \cite{greengard1998an, ferziger1996computational}, and linear elasticity \cite{christiansen1975integral, constanda1995the}. Over the past several decades, numerous traditional numerical methods have been developed to solve equation \eqref{eq:biharmonic}, which can generally be classified into two categories: direct (uncoupled) methods and coupled (decomposed) methods.

For the direct approach, examples include finite difference methods (FDM) based on an uncoupled discretization scheme \cite{lamichhane2014finite, gupta1979direct, altas1998multigrid, ben-artzi2008a, bialecki2012a}, the finite volume method (FVM) \cite{chun-jiabi2004mortar, wang2004a, eymard2012finite}, and the finite element method (FEM), which offers both non-conforming \cite{Baker1977Finite, Lascaux1975Some, Morley2016The} and conforming FEM options \cite{Zienkiewicz2005The, Ciarlet1978The}.

For the indirect approach, which uses a coupled scheme, the splitting method is commonly employed by introducing supplementary attributes and decomposing the biharmonic equation into two Poisson equations. This method further leverages FDM \cite{smith1968the, ehrlich1971solving, bialecki2012a}, FEM, and mixed element techniques \cite{brezzi1991mixed, cheng2000some, davini2000an, lamichhane2011a, stein2019a} to solve the derived equations from the original high-order form. Additionally, collocation methods \cite{mai-duy2009a, bialecki2010spectral, bialecki2020a} and radial basis function (RBF) methods \cite{maiduy2005an, adibi2007numerical, li2011a} are notable approaches for solving \eqref{eq:biharmonic}.

However, these numerical methods are computationally demanding due to the need for assembling and solving large, sparse systems of equations. Meshing and numerical integration in FEM, or grid generation and resolution in FDM, present challenges when balancing accuracy and computational efficiency. Moreover, traditional (or ``classical'') numerical methods, which exclude machine learning-based approaches, often struggle with complexities such as irregular domains and high-dimensional spaces, commonly referred to as the Curse of Dimensionality.

Deep learning, particularly through deep neural networks (DNNs), has shown remarkable success in tackling mathematical problems encountered in scientific computing and engineering, due to its powerful capacity for nonlinear approximation. DNNs have been effectively applied to solve a range of problems, including time-independent and time-dependent PDEs~\cite{e2017deep, e2018the, berg2018a, sirignano2018dgm}, stochastic differential equations \cite{nabian2019a}, inverse problems \cite{raissi2019physics}, and molecular modeling \cite{zou2020deep}. Notably, the Deep Ritz Method (DRM) \cite{e2018the} and Physics-Informed Neural Networks (PINN) \cite{raissi2019physics} have gained significant attention for solving diverse partial differential equations (PDEs), achieving impressive performance in these areas. Utilizing DNNs for PDE approximation provides several advantages, including mesh-free implementation without domain discretization, flexibility in handling complex problems, and the capability for parallel computing to address large-scale tasks \cite{LIU2023272, Dwivedi2020}. Generally, however, PINN demands substantial computational resources for high-order problems to obtain various-order derivatives in the governing PDEs, especially in high-dimensional settings \cite{han2017deep}. Moreover, PINN's performance can degrade under complex or implicit boundary constraints in high-dimensional spaces~\cite{berg2018unified}.

To address these challenges, \citet{lyu2022mim} introduced a deep mixed residual method (MIM) for high-order PDEs by decomposing them into multiple lower-order systems, using a shared DNN to approximate both the solution and its derivatives. Nonetheless, this approach does not inherently align with the original PDE mechanisms, impacting MIM's quality. Alternatively, another MIM formulation splits the fourth-order problem \eqref{eq:biharmonic} into several second-order problems, yielding favorable results. For the biharmonic problem \eqref{eq:biharmonic} with Navier boundary constraints, a coupled deep neural network (CDNN) based on DRM architecture was developed by applying a coupled scheme for \eqref{eq:biharmonic}, configured with two independent Fourier-induced DNNs. However, this configuration requires coercive boundary conditions, making it unsuitable for Dirichlet biharmonic equations \cite{li2022solving}.

The Fourier feature mapping technique, which introduces sine and cosine functions as activation functions before the first hidden layer in a DNN, has been shown to improve the performance of PINN. Due to the innate properties of DNNs, this model is particularly sensitive to low-frequency components of the target function, effectively capturing these initially while having limited capacity for high-frequency components~\cite{rahaman2018spectral, xu2020frequency, wang2020eigenvector, tancik2020fourier, li2023deep}. This mapping approach also preserves the essential properties of the double series expansion method used to solve the biharmonic equation \eqref{eq:biharmonic}. By Fourier-transforming the input in the first layer, the model preserves each input frequency's characteristics through subsequent iterations \cite{Wang2021Eigenvector, tancik2020fourier}. Additionally, Fourier activation functions retain frequency information even after multiple derivatives, unlike traditional activation functions like tanh or sin, making them well-suited for solving multi-scale and high-order problems.

In this paper, we aim to improve the performance of the coupled physics-informed deep neural network (CPINN) method for solving the high-order elliptic equation \eqref{eq:biharmonic} by combining its coupled scheme with Fourier feature mapping, resulting in the Fourier-CPINN (FCPINN) model. Inspired by the coupled scheme, we avoid high-order differentiation by introducing auxiliary functions to decompose the biharmonic equation into two Poisson equations, thereby reducing the problem's complexity. We then integrate the CPINN method with Fourier feature mapping, forming FCPINN. Inspired by double-triangle series and Fourier expansion, this approach provides a robust solution for capturing the complex, multi-frequency characteristics of biharmonic equations. Numerical results reveal that the Fourier transform-based activation function significantly enhances the model's performance, as also observed in \citet{wang2020eigenvector} and \citet{Matthew2020Fourier}, highlighting the benefits of such activation functions. The main contributions of this paper are as follows:

\begin{enumerate}
    \item Utilizing the coupling properties of the fourth-order biharmonic equation and the implicit nature of neural networks, we developed a multi-output FCPINN architecture for solving this equation. This architecture effectively addresses problems with Dirichlet and Navier boundary conditions through the deep Ritz variational method. Compared to existing DNN methods, this architecture produces coupled outputs, reducing computation time compared to multi-network methods. The model achieves significant improvements in accuracy and efficiency over traditional PINN and DRM models.
    
    \item To mitigate the spectral bias of DNNs, we apply Fourier feature mapping in the first layer with sine and cosine activation functions. This approach mimics Fourier expansion, where the first layer acts as Fourier-based functions, combining nonlinearly to enhance the DNN's outputs. Comparative studies indicate that Fourier feature mapping markedly improves FCPINN's adaptability to large-scale input data, addressing multi-scale and high-order challenges in solving biharmonic equations.
    
    \item Through simulations of high-dimensional biharmonic problems, we demonstrate that our FCPINN model achieves considerably higher accuracy while requiring fewer computational resources across various dimensional spaces, such as elasticity mechanics, as shown in the numerical experiments.
\end{enumerate}

The paper is organized as follows. Section \ref{sec:02} provides a brief overview of the DNN, PINN frameworks, and their formulations. Section \ref{sec:03} details the construction of the FCPINN architecture, based on the coupled scheme, to approximate the biharmonic equation's solution with Fourier mapping as the activation function. Section \ref{sec:04} presents numerical experiments that assess the FCPINN model's performance across different boundary conditions and high-dimensional, high-order problems. Finally, Section \ref{sec:05} offers a brief conclusion.

\section{Preliminaries}\label{sec:02}
In this section, we provide a detailed overview of the mathematical concepts and formulas relevant to Deep Neural Networks (DNN) and Physics-Informed Neural Networks (PINN). 

\subsection{Deep Neural Networks}\label{sec:dnn}
We begin by introducing the foundational concepts and formulations of Deep Neural Networks (DNN) to give readers a functional understanding of DNN structures. Mathematically, a DNN establishes a mapping relationship between nodes across layers:
\begin{equation}
\mathcal{F}: \bm{x}\in\mathbb{R}^{d}\Longrightarrow \bm{y}=\mathcal{F}(x)\in\mathbb{R}^{c}
\end{equation}
where $d$ and $c$ denote the dimensions of the input and output, respectively. A standard neuron unit in a DNN defines the following mapping relationship between its input $\bm{x} \in \mathbb{R}^{d}$ and output $\bm{y} \in \mathbb{R}^{m}$ as shown in~\eqref{eq0201}:
\begin{equation}\label{eq0201}
\bm{y}=\{y_1, y_2, \cdots, y_m\}~~\textup{and}~~  y_l= \sigma\left(\sum_{n=1}^{d} w_{ln}*x_n  + b_l\right)
\end{equation}
In this setup, $y_1$ represents the first layer of the neural network, comprising several nodes. The weight and bias for the $l^{\text{th}}$ neuron in this layer are denoted by $w_{ln}$ and $b_l$, respectively. The function $\sigma$ is a nonlinear operator, also known as the activation function, that transforms each element. With the first layer defined as above, a deep neural network (DNN) is constructed through layer-by-layer transformations of each node, utilizing linear calculations (weights and biases) alongside nonlinear activation functions. For an input $\bm{x} \in \mathbb{R}^d$, the mathematical form of this mapping relationship is given by:
\begin{equation}\label{form2dnn}
	\begin{cases}
		\bm{y}^{[0]} = \bm{x}\\
		\bm{y}^{[\ell]} = \sigma\circ(\bm{W}^{[\ell]}\bm{y}^{[\ell-1]}+\bm{b}^{[\ell]}), ~~\text{for}~~\ell =1, 2, 3, \cdots\cdots, L
	\end{cases}.
\end{equation}
Here, $\bm{W}^{[\ell]} \in \mathbb{R}^{n_{\ell+1}\times n_{\ell}}$ and $\bm{b}^{[\ell]}\in\mathbb{R}^{n_{\ell+1}}$ represent the weights and biases of the $\ell$-th hidden layer. The transformation relationship between hidden layers is similar to that of the first layer. Using the previously hidden layer as input, each node's values in the next layer are computed through operations involving linear weights, biases, and nonlinear activation functions. The notation ``$\circ$'' denotes an element-wise operation for each node. Information transfer between hidden layers is modulated by each node's weights and biases, transformed by the activation function, and passed to corresponding nodes in the next layer. This iterative process defines the structure of the deep neural network.

The collection of weights and biases for all nodes forms the neural network’s parameter set: $\bm{W}^{[1]},\cdots, \bm{W}^{[L]}$, $\bm{b}^{[1]},\cdots, \bm{b}^{[L]}$, which we denote by $\bm{\theta}$. The output of the DNN, denoted as $\bm{y}(\bm{x};\bm{\theta})$, represents the network's output set.

\begin{remark}
In general, $\bm{\theta}$, which comprises all weights and biases in the DNN, is initialized using the Xavier initialization method as described in \citet{kumar2017weight}, i.e., from the following normal distribution:
\begin{equation*}
    \mathcal{N}\left(0,\left(\frac{2}{n_{in}+n_{out}}\right)^2\right),
\end{equation*}
where $n_{in}$ and $n_{out}$ denote the dimensions of the input and output of the corresponding layer, respectively.
\end{remark}

\subsection{Physics-Informed Neural Networks}\label{sec:0202}

Consider, for instance, a set of parameterized PDEs expressed as:
\begin{equation}\label{eq2PPDE}
\begin{aligned}
&\mathcal{N}_{\bm{\lambda}}[\hat{u}(\bm{x})]=\hat{f}(\bm{x}), \quad \bm{x} \in \Omega \\
&\mathcal{B}\hat{u}\left(\bm{x}\right)=\hat{g}(\bm{x}), \quad \bm{x} \in \partial \Omega\\
\end{aligned}
\end{equation}
where $\mathcal{N}_{\bm{\lambda}}$ is a partial differential operator with parameters $\bm{\lambda}$ that may represent both linear and nonlinear constraints, and $\mathcal{B}$ is the relevant boundary operator. Here, $\Omega$ and $\partial \Omega$ denote the domain and its boundary, respectively. DNN models are often utilized in PINNs to solve PDEs \eqref{eq2PPDE} by minimizing the following loss function to obtain an optimal solution to the associated PDEs:
\begin{equation}\label{loss2PINN}
Loss=Loss_{PDE}+\gamma Loss_{BC}
\end{equation}
with
\begin{equation}\label{subloss2PINN}
\begin{aligned}
&Loss_{PDE} = \frac{1}{N_P}\sum_{i=1}^{N_P}\left| \mathcal{N}_{\bm{\lambda}}[\hat{u}_{NN}(\bm{x}^i)]-\hat{f}(\bm{x}^i)\right|^2 \\
&Loss_{BC} = \frac{1}{N_B}\sum_{i=1}^{N_B}\bigg{|}\mathcal{B}\hat{u}_{NN}\left(\bm{x}^i\right)-\hat{g}(\bm{x}^i)\bigg{|}^2
\end{aligned}
\end{equation}
where $\gamma$ is a weighting parameter for the boundary loss. Here, $Loss_{PDE}$ and $Loss_{BC}$ denote the residual errors for the governing equations and the specified boundary conditions (BC), respectively. If additional data are available within the domain, an extra loss term can be added to account for the difference between the approximation and the raw data:
\begin{equation}
    Loss_{D} = \frac{1}{N_D}\sum_{i=1}^{N_D}\bigg{|}\hat{u}_{NN}(\bm{x}^i)-u_{Data}^i\bigg{|}^2.
\end{equation}

\begin{figure}[!htbp]
    \centering
    \includegraphics[scale=0.4]{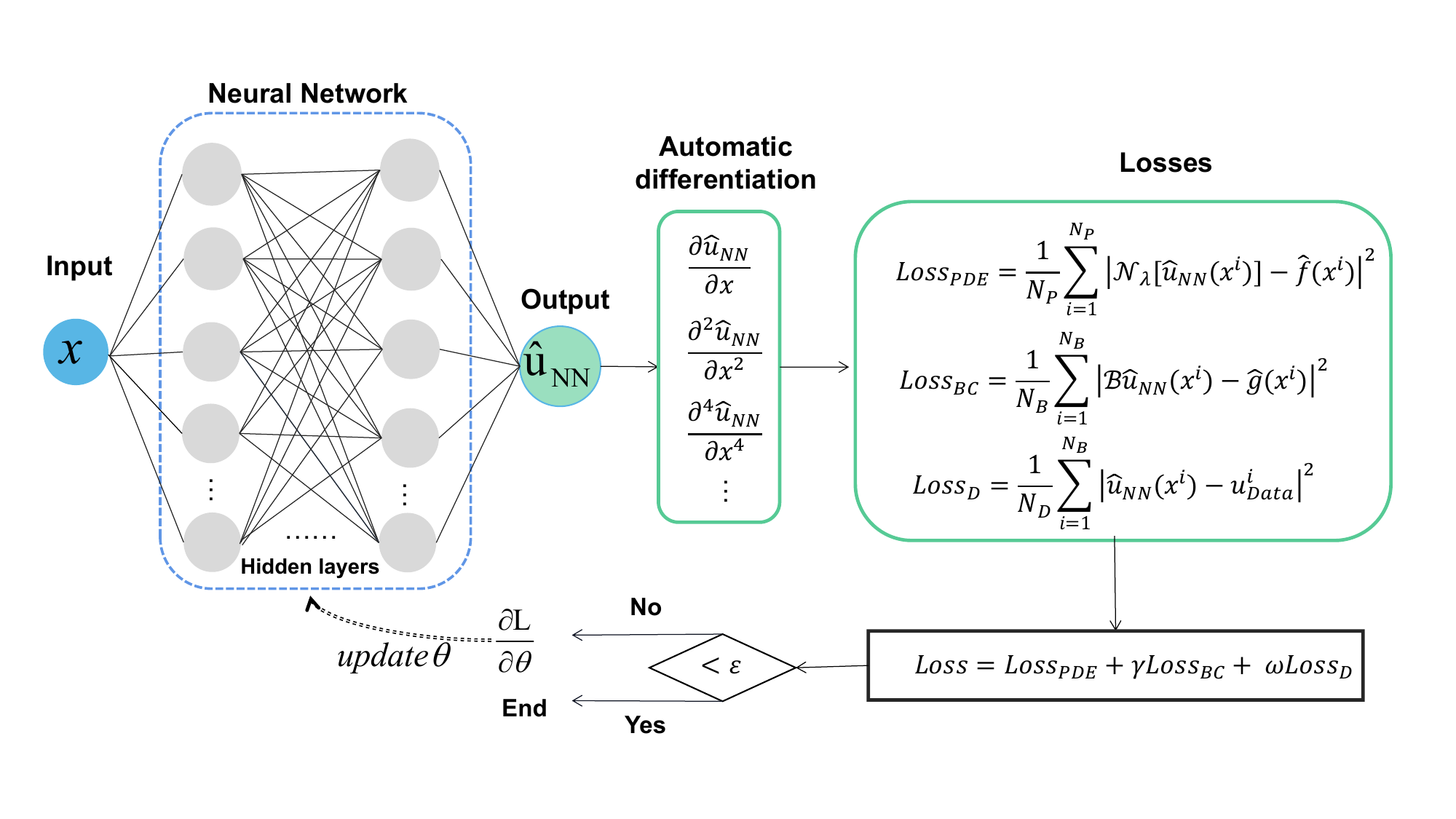}
    \caption {Illustration of the PINN structure.}
    \label{fig: PINN Structure}
\end{figure}

In solving parameterized PDEs, the structure of PINN for \eqref{eq2PPDE} is illustrated in Figure~\ref{fig: PINN Structure}.

\section{Unified Architecture of FCPINN for Biharmonic Equations}\label{sec:03}
\subsection{FCPINN Algorithm for Biharmonic Equations}\label{algor2CELM}

We now introduce the specific approach for using Coupled PINN (CPINN) to solve the biharmonic equation. This method draws on decoupling techniques from traditional numerical methods, such as the Finite Element Method (FEM) and Finite Difference Method (FDM). When handling higher-order problems, we introduce an auxiliary variable $v = \Delta u$, allowing us to rewrite the original fourth-order equation as two second-order Poisson equations, thereby simplifying the computation. These two Poisson equations are expressed as follows:

\begin{equation}\label{couple2Dirichlet}
   \begin{cases}
    \Delta v(\bm{x}) = f(\bm{x},u,\Delta u), ~~\text{in}~~ \Omega\\
     ~--,   ~~~~ \text{on}~~ \partial \Omega
   \end{cases}
   ~~~~\text{and}~~~~
   \begin{cases}
   \Delta u(\bm{x})=v(\bm{x}), ~~\text{in}~~ \Omega\\
   ~~~u(\bm{x})=g(\bm{x}),~~~ \text{on}~~ \partial \Omega\\
   ~~~\frac{\partial u(\bm{x})}{\partial \vec{n}}= h(\bm{x}), ~~\text{on}~ \partial\Omega
   \end{cases}.
\end{equation}
for Dirichlet boundary conditions, where $k(x)$ specifies $v(x)$ on the boundary, and $h(x)$ defines the normal derivative of $u(x)$ on the boundary.

\begin{equation}\label{couple2Navier}
   \begin{cases}
    \Delta v(\bm{x}) = f(\bm{x},u,\Delta u), ~~\text{in}~~ \Omega\\
    ~~~v(\bm{x})=k(\bm{x}),   ~~~~ \text{on}~~ \partial \Omega
   \end{cases}
   ~~~~\text{and}~~~~
   \begin{cases}
   \Delta u(\bm{x})=v(\bm{x}), ~~\text{in}~~ \Omega\\
   ~~~u(\bm{x})=g(\bm{x}),~~~ \text{on}~~ \partial \Omega
   \end{cases}.
\end{equation}
for Navier boundary conditions.

Rather than seeking a single solution to the original problem \eqref{eq:biharmonic}, our goal is to approximate the pair $(u, v)$ that satisfies specific boundary or constraint conditions, minimizing computational load and complexity. Numerical solvers can be used to resolve \eqref{couple2Dirichlet} or \eqref{couple2Navier} by identifying values within the domain $\Omega$ that minimize the least-squares value of the following expression:
\begin{equation}\label{loss2continous}
u^*,v^* = \underset{(u,v)\in \mathbb{H}_0^2(\Omega)\times \mathbb{H}_0^2(\Omega)}{\arg\min}\mathcal{L}(u,v)
\end{equation}
with
\begin{equation}\label{intergalForm}
\mathcal{L}(u,v) = \int_{\Omega}\big{|}\Delta v(\bm{x}) - f(\bm{x},u,\Delta u)\big{|}^2d\bm{x} + \int_{\Omega}\big{|}\Delta u(\bm{x})-v(\bm{x})\big{|}^2d\bm{x}
\end{equation}
where $u, v \in \mathbb{H}_0^2(\Omega)$ are trial functions, and $\mathbb{H}_0^2(\Omega)$ is the admissible function space based on the differential operator $\Delta$ and the boundary operator $\mathcal{B}$. 

Typically, two separate neural networks are used to handle the auxiliary variables $v$ and $u$. However, for equations with Dirichlet boundary constraints, $v$ serves as an intermediate variable without explicit constraint conditions. Given the strong linear and nonlinear approximation abilities of DNNs, they can be applied to solve various linear and nonlinear complex functions. In this study, we employ a single DNN with multiple outputs to approximate $v$ and $u$, denoted as $v_{NN}$ and $u_{NN}$.

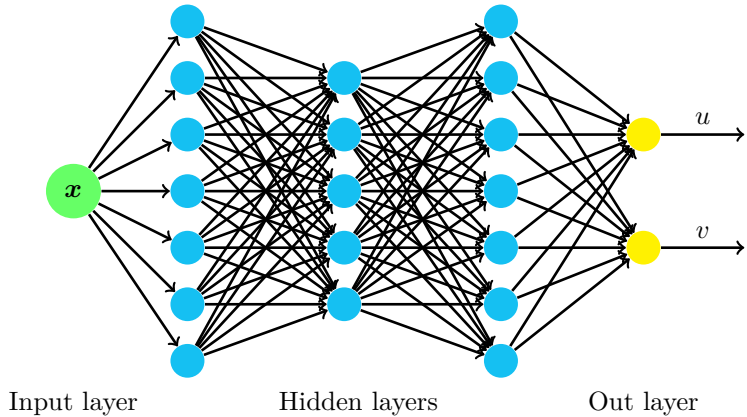
\begin{figure}[H]
    \centering
    \begin{tikzpicture}[scale=0.75]
    \node[] (input) at (2, -3.75) {Input layer};
    
    \node[circle, fill=green!60,inner sep=4.5pt] (x) at (2.0, 0) {$\bm{x}$};		
    
    \node[circle, fill=cyan!70,inner sep=4.5pt] (h10) at (4, 3) {};
    \node[circle, fill=cyan!70,inner sep=4.5pt] (h11) at (4, 2) {};
    \node[circle, fill=cyan!70,inner sep=4.5pt] (h12) at (4, 1) {};
    \node[circle, fill=cyan!70,inner sep=4.5pt] (h13) at (4, 0) {};
    \node[circle, fill=cyan!70,inner sep=4.5pt] (h14) at (4, -1) {};
    \node[circle, fill=cyan!70,inner sep=4.5pt] (h15) at (4, -2) {};
    \node[circle, fill=cyan!70,inner sep=4.5pt] (h16) at (4, -3) {};

    \draw[line width=1.0pt,->] (x) -- (h10);
    \draw[line width=1.0pt,->] (x) -- (h11);
    \draw[line width=1.0pt,->] (x) -- (h12);
    \draw[line width=1.0pt,->] (x) -- (h13);
    \draw[line width=1.0pt,->] (x) -- (h14);
    \draw[line width=1.0pt,->] (x) -- (h15);
    \draw[line width=1.0pt,->] (x) -- (h16);

    \node[circle, fill=cyan!70,inner sep=4.5pt] (h20) at (6.75, 2) {};
    \node[circle, fill=cyan!70,inner sep=4.5pt] (h21) at (6.75, 1) {};
    \node[circle, fill=cyan!70,inner sep=4.5pt] (h22) at (6.75, 0) {};
    \node[circle, fill=cyan!70,inner sep=4.5pt] (h23) at (6.75, -1) {};
    \node[circle, fill=cyan!70,inner sep=4.5pt] (h24) at (6.75, -2) {};
    \node[] (input) at (7, -3.75) {Hidden layers};
    
    \draw[line width=1.0pt,->] (h10) -- (h20);
    \draw[line width=1.0pt,->] (h10) -- (h21);
    \draw[line width=1.0pt,->] (h10) -- (h22);
    \draw[line width=1.0pt,->] (h10) -- (h23);
    \draw[line width=1.0pt,->] (h10) -- (h24);
    
    \draw[line width=1.0pt,->] (h11) -- (h20);
    \draw[line width=1.0pt,->] (h11) -- (h21);
    \draw[line width=1.0pt,->] (h11) -- (h22);
    \draw[line width=1.0pt,->] (h11) -- (h23);
    \draw[line width=1.0pt,->] (h11) -- (h24);
    
    \draw[line width=1.0pt,->] (h12) -- (h20);
    \draw[line width=1.0pt,->] (h12) -- (h21);
    \draw[line width=1.0pt,->] (h12) -- (h22);
    \draw[line width=1.0pt,->] (h12) -- (h23);
    \draw[line width=1.0pt,->] (h12) -- (h24);
    
    \draw[line width=1.0pt,->] (h13) -- (h20);
    \draw[line width=1.0pt,->] (h13) -- (h21);
    \draw[line width=1.0pt,->] (h13) -- (h22);
    \draw[line width=1.0pt,->] (h13) -- (h23);
    \draw[line width=1.0pt,->] (h13) -- (h24);;
    
    \draw[line width=1.0pt,->] (h14) -- (h20);
    \draw[line width=1.0pt,->] (h14) -- (h21);
    \draw[line width=1.0pt,->] (h14) -- (h22);
    \draw[line width=1.0pt,->] (h14) -- (h23);
    \draw[line width=1.0pt,->] (h14) -- (h24);
    
    \draw[line width=1.0pt,->] (h15) -- (h20);
    \draw[line width=1.0pt,->] (h15) -- (h21);
    \draw[line width=1.0pt,->] (h15) -- (h22);
    \draw[line width=1.0pt,->] (h15) -- (h23);
    \draw[line width=1.0pt,->] (h15) -- (h24);
    
    \draw[line width=1.0pt,->] (h16) -- (h20);
    \draw[line width=1.0pt,->] (h16) -- (h21);
    \draw[line width=1.0pt,->] (h16) -- (h22);
    \draw[line width=1.0pt,->] (h16) -- (h23);
    \draw[line width=1.0pt,->] (h16) -- (h24);
    
    \node[circle, fill=cyan!70,inner sep=4.5pt] (h30) at (9.5, 3) {};
    \node[circle, fill=cyan!70,inner sep=4.5pt] (h31) at (9.5, 2) {};
    \node[circle, fill=cyan!70,inner sep=4.5pt] (h32) at (9.5, 1) {};
    \node[circle, fill=cyan!70,inner sep=4.5pt] (h33) at (9.5, 0) {};
    \node[circle, fill=cyan!70,inner sep=4.5pt] (h34) at (9.5, -1) {};
    \node[circle, fill=cyan!70,inner sep=4.5pt] (h35) at (9.5, -2) {};
    \node[circle, fill=cyan!70,inner sep=4.5pt] (h36) at (9.5, -3) {};	
    
    \draw[line width=1.0pt,->] (h20) -- (h30);
    \draw[line width=1.0pt,->] (h20) -- (h31);
    \draw[line width=1.0pt,->] (h20) -- (h32);
    \draw[line width=1.0pt,->] (h20) -- (h33);
    \draw[line width=1.0pt,->] (h20) -- (h34);
    \draw[line width=1.0pt,->] (h20) -- (h35);
    \draw[line width=1.0pt,->] (h20) -- (h36);
    
    \draw[line width=1.0pt,->] (h21) -- (h30);
    \draw[line width=1.0pt,->] (h21) -- (h31);
    \draw[line width=1.0pt,->] (h21) -- (h32);
    \draw[line width=1.0pt,->] (h21) -- (h33);
    \draw[line width=1.0pt,->] (h21) -- (h34);
    \draw[line width=1.0pt,->] (h21) -- (h35);
    \draw[line width=1.0pt,->] (h21) -- (h36);
    
    \draw[line width=1.0pt,->] (h22) -- (h30);
    \draw[line width=1.0pt,->] (h22) -- (h31);
    \draw[line width=1.0pt,->] (h22) -- (h32);
    \draw[line width=1.0pt,->] (h22) -- (h33);
    \draw[line width=1.0pt,->] (h22) -- (h34);
    \draw[line width=1.0pt,->] (h22) -- (h35);
    \draw[line width=1.0pt,->] (h22) -- (h36);
    
    \draw[line width=1.0pt,->] (h23) -- (h30);
    \draw[line width=1.0pt,->] (h23) -- (h31);
    \draw[line width=1.0pt,->] (h23) -- (h32);
    \draw[line width=1.0pt,->] (h23) -- (h33);
    \draw[line width=1.0pt,->] (h23) -- (h34);
    \draw[line width=1.0pt,->] (h23) -- (h35);
    \draw[line width=1.0pt,->] (h23) -- (h36);
    
    \draw[line width=1.0pt,->] (h24) -- (h30);
    \draw[line width=1.0pt,->] (h24) -- (h31);
    \draw[line width=1.0pt,->] (h24) -- (h32);
    \draw[line width=1.0pt,->] (h24) -- (h33);
    \draw[line width=1.0pt,->] (h24) -- (h34);
    \draw[line width=1.0pt,->] (h24) -- (h35);
    \draw[line width=1.0pt,->] (h24) -- (h36);
    
    \node[circle, fill=yellow!100,inner sep=4.5pt] (o0) at (12, 1.0) {};	
    \node[circle, fill=yellow!100,inner sep=4.5pt] (o1) at (12, -1.0) {};	
    
    \draw[line width=1.0pt,->] (h30) -- (o0);
    \draw[line width=1.0pt,->] (h30) -- (o1);
    
    \draw[line width=1.0pt,->] (h31) -- (o0);
    \draw[line width=1.0pt,->] (h31) -- (o1);
    
    \draw[line width=1.0pt,->] (h32) -- (o0);
    \draw[line width=1.0pt,->] (h32) -- (o1);
    
    \draw[line width=1.0pt,->] (h33) -- (o0);
    \draw[line width=1.0pt,->] (h33) -- (o1);
    
    \draw[line width=1.0pt,->] (h34) -- (o0);
    \draw[line width=1.0pt,->] (h34) -- (o1);
    
    \draw[line width=1.0pt,->] (h35) -- (o0);
    \draw[line width=1.0pt,->] (h35) -- (o1);
    
    \draw[line width=1.0pt,->] (h36) -- (o0);
    \draw[line width=1.0pt,->] (h36) -- (o1);
    
    \node[circle] (y0) at (14, 1) {};
    \node[circle] (y1) at (14, -1) {};
    \node[] (input) at (12, -3.75) {Out layer};
    
    \draw[line width=1.0pt,->] (o0) -- node[above]{$u$}(y0);
    \draw[line width=1.0pt,->] (o1) -- node[above]{$v$}(y1);
    \end{tikzpicture}
    \caption{ \small The DNN for predicting the state and flux parameters  with multi-output}
    \label{fig2multioutput}
\end{figure}

In Figure~\ref{fig2multioutput}, we illustrate a three-layer neural network with multi-output, using $\bm{x} \in \mathbb{R}^2$ as an example input.

For the established expressions of the auxiliary variable $v$ and the target function $u$, the Monte Carlo method \cite{robert1999monte} is commonly employed to discretize \eqref{intergalForm}. We then apply the PINN method to carry out the computations, leading to the following loss function:
\begin{equation}\label{loss2ellptic}
\begin{aligned}
    \mathcal{L}_{in}(\mathcal{S}_I;\bm{\theta})&=\frac{|\Omega|}{N_{in}}\sum_{i=1}^{N_{in}}\Big{|}\Delta v_{NN}(\bm{x}_I^i;\bm{\theta}) - f(\bm{x}_I^i,u_{NN}(\bm{x}_I^i;\bm{\theta}),\Delta u_{NN}(\bm{x}_I^i;\bm{\theta}))\Big{|}^2\\
    &+\frac{|\Omega|}{N_{in}}\sum_{i=1}^{N_{in}}\Big{|}\Delta u_{NN}(\bm{x}_I^i;\bm{\theta}) - v_{NN}(\bm{x}_I^i;\bm{\theta})\Big{|}^2 ~~\text{for}~~\bm{x}_I^i \in S_I.
\end{aligned}
\end{equation}

In this context and henceforth, $S_I$ represents the sampled points from $\Omega$ with a specified probability density.

Boundary conditions are essential constraints for numerical methods like the Finite Difference Method (FDM) and Finite Element Method (FEM) in solving PDE problems, ensuring solution uniqueness and accuracy. Similarly, enforcing boundary conditions within the DNN representation is a crucial aspect. For the Dirichlet boundary, we know that $\partial u/\partial \vec{n}=\nabla u\cdot \vec{n}$, where $\vec{n}$ represents the unit outward normal on the boundary of $\Omega$. Under the boundary constraints of \eqref{couple2Dirichlet} and \eqref{couple2Navier}, the expressions of $v_{NN}$ and $u_{NN}$ on the boundary of $\Omega$ should satisfy
\begin{equation}\label{loss_bd2dirichlet}
\mathcal{L}_{bd}(\mathcal{S}_B;\bm{\theta})=\frac{1}{M}\sum_{j=1}^{M}\Big{|} v_{NN}(\bm{x}_B^j;\bm{\theta}) - g(\bm{x}_B^j)\Big{|}^2+\frac{1}{M}\sum_{j=1}^{M} \Big{|} \nabla u_{NN}(\bm{x}_B^h;\bm{\theta})\cdot\vec{n} - h(\bm{x}_B^j)\Big{|}^2 ~~\text{for}~~\bm{x}_B^j\in \partial \Omega
\end{equation}
for Dirichlet boundary conditions, or
\begin{equation}\label{loss_bd2navier}
\mathcal{L}_{bd}(\mathcal{S}_B;\bm{\theta})=\frac{1}{M}\sum_{j=1}^{M}\Big{|} v_{NN}(\bm{x}_B^j;\bm{\theta}) - g(\bm{x}_B^j)\Big{|}^2+\frac{1}{M}\sum_{j=1}^{M} \Big{|}  v_{NN}(\bm{x}_B^h;\bm{\theta}) - k(\bm{x}_B^j)\Big{|}^2 ~~\text{for}~~\bm{x}_B^j\in \partial \Omega
\end{equation}
for Navier boundary conditions. Here, $S_B$ represents the sampled points on $\partial \Omega$ with a specified probability density.

To optimize the parameters in the DNN, the following cost function should be minimized:
\begin{equation*}
	\mathcal{L}(\left\{\mathcal{S}_I,\mathcal{S}_B\right\};\bm{\theta}) = \mathcal{L}_{in}(\mathcal{S}_I;\bm{\theta}) + \gamma\mathcal{L}_{bd}(\mathcal{S}_B;\bm{\theta}).
\end{equation*}
In this context, ${\mathcal{S}_I}=\{\bm{x}_I^i\}_{i=1}^{N}$ and ${\mathcal{S}_B}=\{\bm{x}_B^j\}_{j=1}^{M}$ denote the training data within the internal region $\Omega$ and on the boundary $\partial \Omega$, respectively. The penalty parameter $\gamma$ adjusts the weight of the loss contribution in the boundary region $\partial \Omega$, gradually increasing during iteration to improve boundary fitting.

If $\mathcal{L}(\bm{\theta})$ approaches a sufficiently small value, $L(\cdot,\bm{\theta})$ will approximate the solution of the following equation \eqref{eq:biharmonic}. Thus, we minimize $\mathcal{L}(\bm{\theta})$ via the approximation function $L(\cdot,\bm{\theta})$ to find a set of parameters $\bm{\theta}$ that satisfy the conditions:
\begin{equation*}
\bm{\theta}^* = \arg\min \mathcal{L}(\cdot;\bm{\theta})~~\Longleftrightarrow~~ u(\bm{x})\approx u_{NN}(\bm{x};\bm{\theta}^*).
\end{equation*}

To obtain the parameter set $\bm{\theta}^*$, we can use gradient descent (GD) to update $\bm{\theta}$. This method can be applied to all training samples simultaneously or to a subset of samples in each iteration. Stochastic gradient descent (SGD) and its improvements, such as Adam~\cite{kingma2015adam}, RMSprop~\cite{hinton2012neural}, and Adagrad~\cite{duchi2011adaptive}, are commonly used optimization techniques in deep learning. Unlike conventional GD, SGD evaluates only $n$ functions per iteration or selects a small batch for updates rather than iterating over a single item. The parameter update in standard SGD can be expressed as:
\begin{equation*}
\bm{\theta}^{k+1}=\bm{\theta}^{k}-\alpha^k\nabla_{\bm{\theta}^k}\mathcal{L}(\bm{x};\bm{\theta}^{k})~~\text{with}~~\bm{x}\in \mathcal{S}_I~\text{or}~\bm{x}\in\mathcal{S}_B
\end{equation*}
where $\alpha^k$ is the learning rate.

\subsection{Choice of Activation Function for PINN}

\begin{figure}[H]
    \centering
    \subfigure[Gaussian function]{
        \label{Gaussian} 
        \includegraphics[scale=0.33]{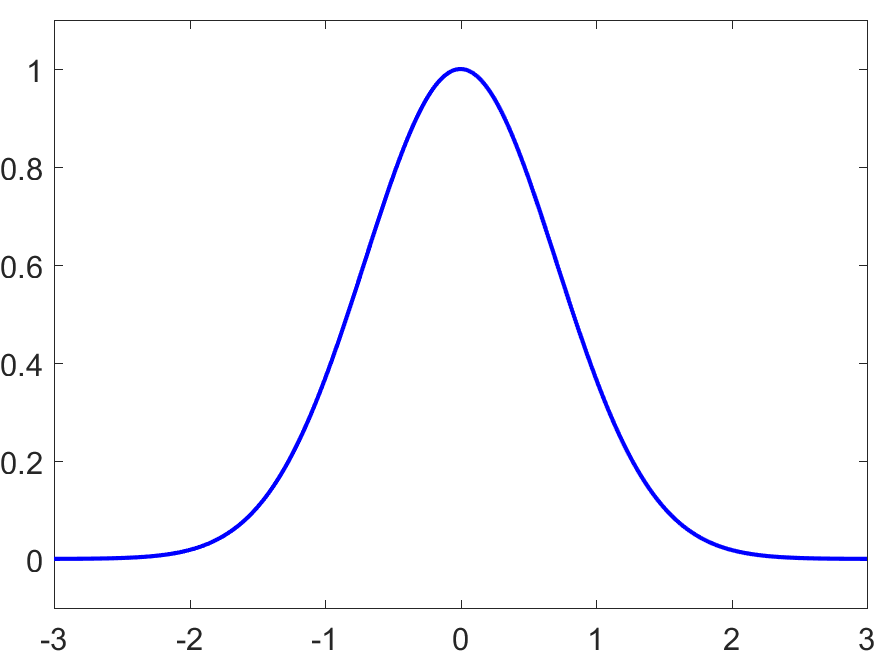}
    }
    \subfigure[1st-order derivative of Gaussian]{
        \label{1stDeri2Gaussian}
        \includegraphics[scale=0.33]{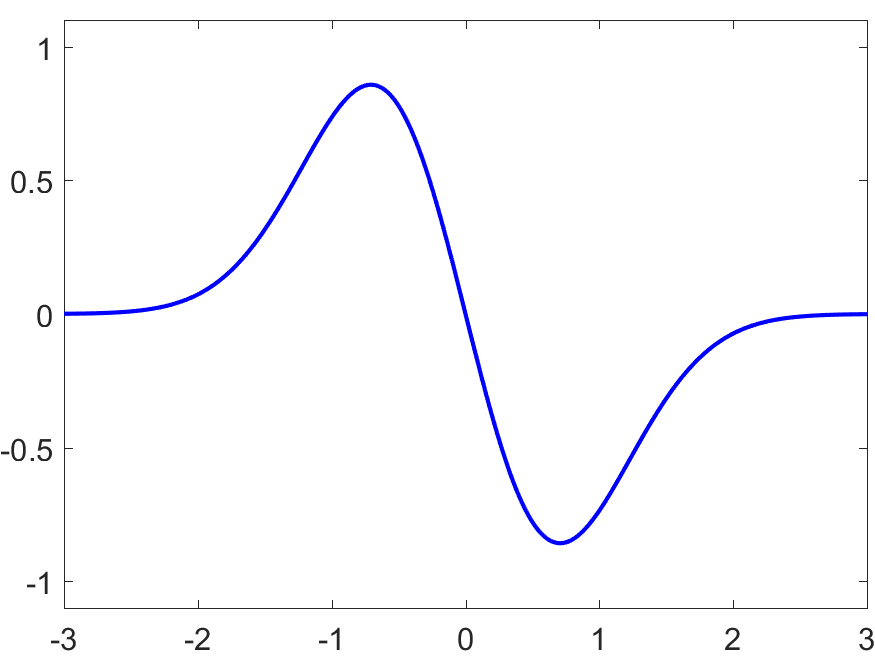}
    }
    \subfigure[2nd-order derivative of Gaussian]{
    \label{2ndDeri2Gaussian}
    \includegraphics[scale=0.33]{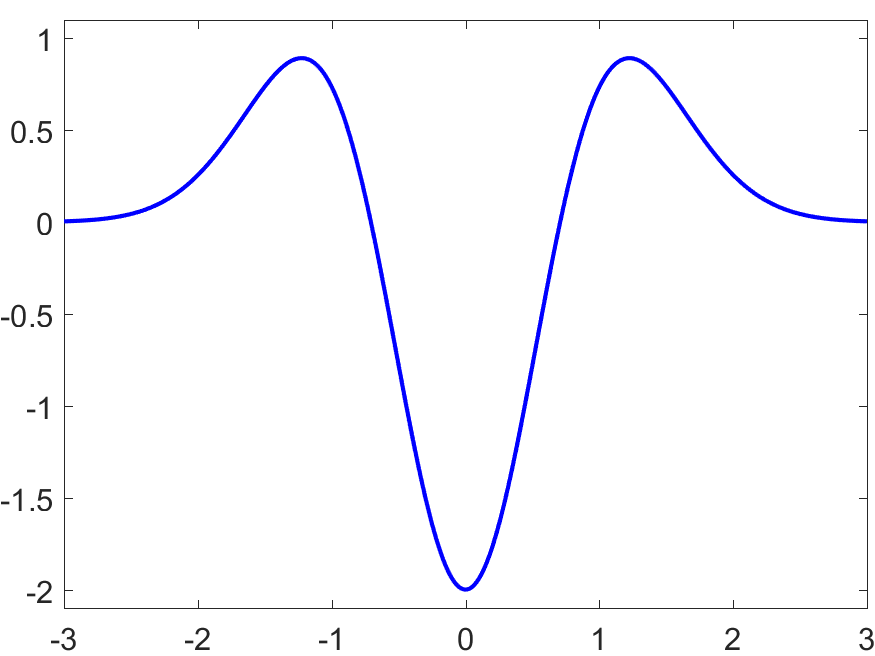}
    }
    \subfigure[4th-order derivative of Gaussian]{
    \label{4thDeri2Gaussian}
    \includegraphics[scale=0.33]{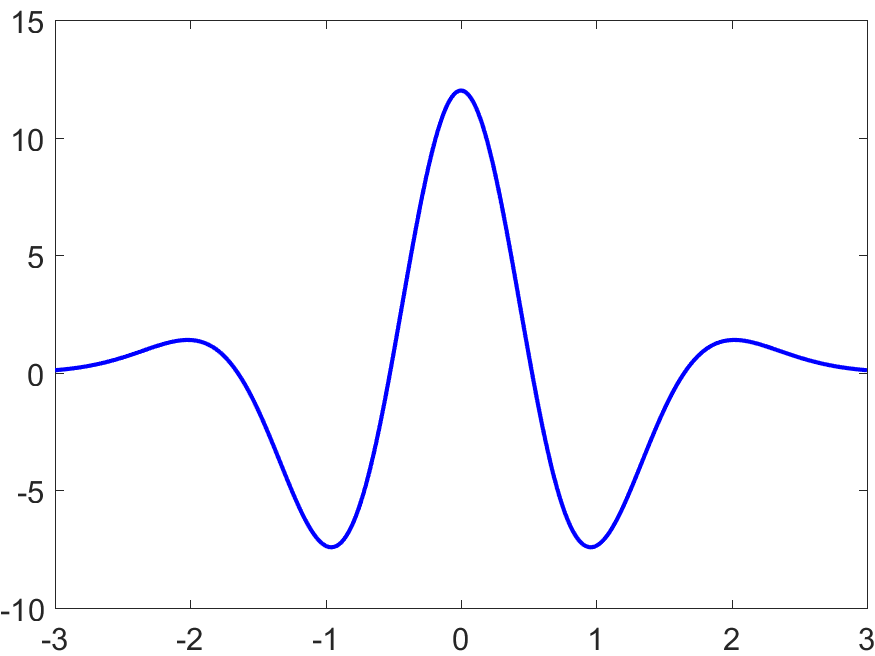}
    }
    \caption{\small The Gaussian function and its 1st-order, 2nd-order, and 4th-order derivatives.}
    \label{gauss_act_func}
\end{figure}

\begin{figure}[H]
    \centering
    \subfigure[Tanh function]{
        \label{Tanh} 
        \includegraphics[scale=0.33]{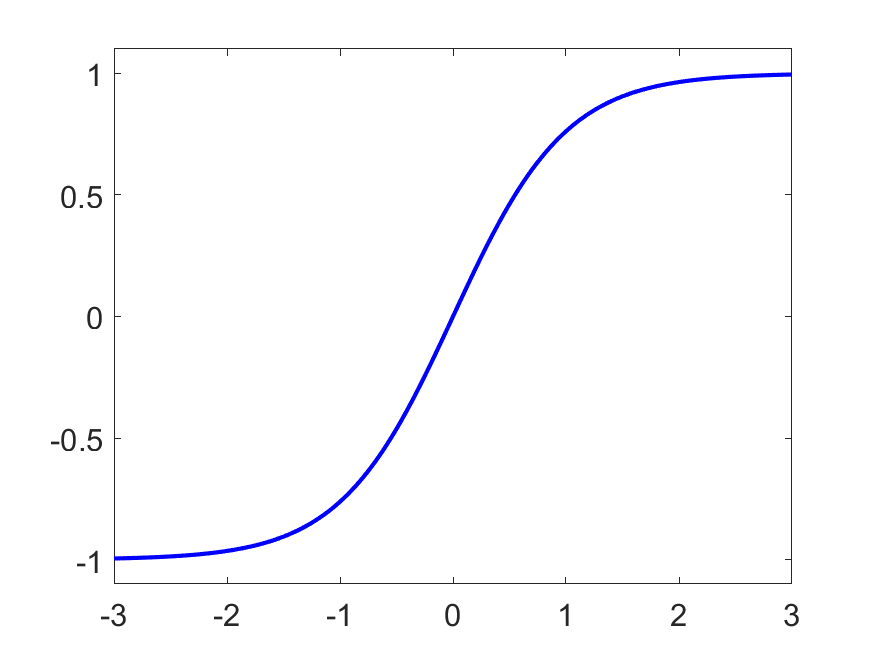}
    }
    \subfigure[1st-order derivative of tanh]{
        \label{1stDeri2tanh}
        \includegraphics[scale=0.33]{figure2avationct/1st_deri2gauss.eps}
    }
    \subfigure[2nd-order derivative of tanh]{
    \label{2ndDeri2tanh}
    \includegraphics[scale=0.33]{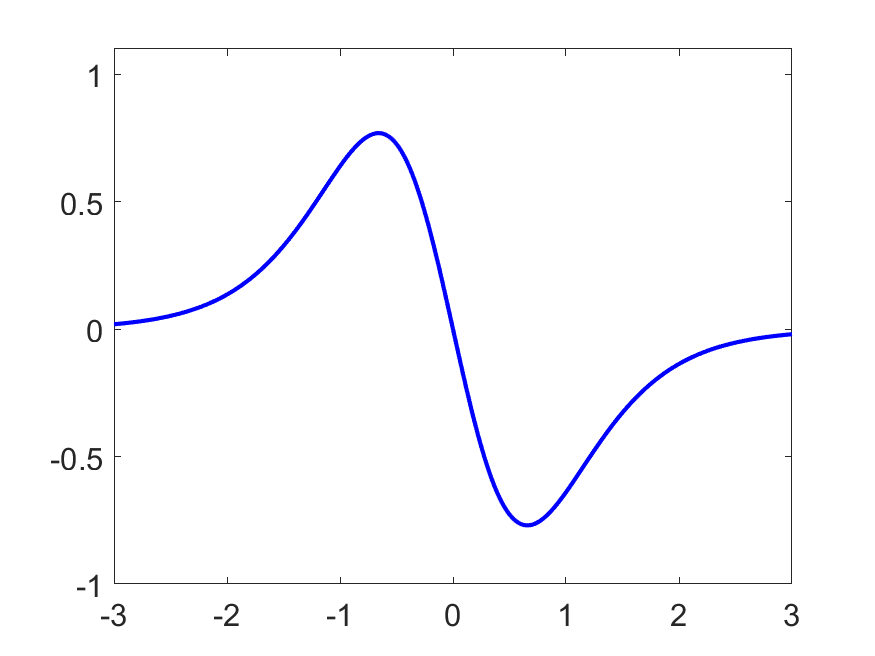}
    }
    \subfigure[4th-order derivative of tanh]{
    \label{4thDeri2tanh}
    \includegraphics[scale=0.33]{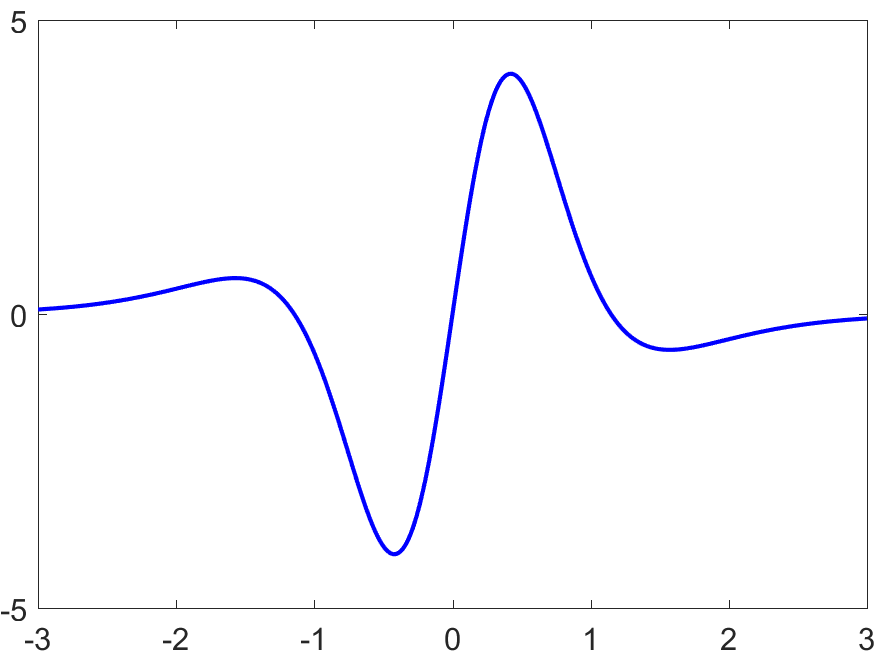}
    }
    \caption{\small The tanh function and its 1st-order, 2nd-order, and 4th-order derivatives.}
    \label{tanh_act_func}
\end{figure}

The design of the activation function is critical to the performance of DNNs. Since the activation function directly influences both the input and internal data flow, its nonlinear characteristics significantly impact the DNN's final performance. Typically, DNNs utilize various activation functions for nonlinear transformations at the input and internal layers. Common activation functions include $\sin(\bm{z})$, $\text{ReLU}(\bm{z}) = \max\{0,\bm{z}\}$, and $\tanh(\bm{z})$. The Gaussian and tanh functions satisfy specific conditions and exhibit strong regularity, ensuring the boundedness of the output in PINN models and enhancing performance. The first, second, and fourth derivatives of the Gaussian function $\sigma(x)=e^{-x^2}$ and the hyperbolic tangent function $\tanh(x) = \frac{e^x - e^{-x}}{e^x + e^{-x}}$ are given by 
\begin{equation}\label{multi_deri2gaussian}
\left\{
\begin{aligned}
\Big{(}\sigma(x)\Big{)}' &= -2xe^{-x^2}\\
\Big{(}\sigma(x)\Big{)}''&= -2e^{-x^2}+4x^2e^{-x^2}\\
\Big{(}\sigma(x)\Big{)}''''&= 12e^{-x^2}-48x^2e^{-x^2}+16x^4e^{-x^2}
\end{aligned}\right.,
\end{equation}
and
\begin{equation}\label{multi_deri2tan}
\left\{
\begin{aligned}
\tanh'(x) &= 1-\tanh^2(x)\\
\tanh''(x) &= -2\tanh(x)+2\tanh^3(x)\\
\tanh''''(x) &= \Big{(}16\tanh(x)-24\tanh^3(x)\Big{)}\Big{(}1-\tanh^2(x)\Big{)}
\end{aligned}\right.,
\end{equation}
respectively. Figures~\ref{gauss_act_func} and~\ref{tanh_act_func} show the curves of the Gaussian and tanh functions along with their respective first, second, and fourth derivatives.

In the above function graphs, for inputs within a narrow range (e.g., within 0-1), the derivatives of the sigmoid, Gaussian, and tanh functions are bounded and converge. However, due to the simplicity of these functions, they may not fully capture information outside this input range. From an approximation standpoint, each node in a DNN's hidden layer contains basic parameters, which are processed into an output by the activation function. The DNN's output is thus a linear or nonlinear combination of these nodes modulated by the activation function. Therefore, designing better activation functions and neural network frameworks remains an open question in enhancing DNN capacity. Standard DNNs are adept at solving low-frequency problems, yet they face challenges when addressing high-frequency oscillations.

Recent research \cite{xu2020frequency,rahaman2018spectral} has identified a phenomenon called spectral bias, or frequency preference, in DNNs. DNNs tend to fit functions by scanning across various frequencies, displaying different sensitivities to information at different frequencies. The neural tangent kernel (NTK) has been used to explain this behavior \cite{jacot2018neural, wang2020eigenvector}. To address spectral bias, Fourier feature mapping can be applied by combining sine and cosine functions as the activation function, enhancing the DNN's capacity to fit high-frequency information by leveraging Fourier transform characteristics \cite{Matthew2020Fourier,wang2020eigenvector}. Additionally, using a global Fourier mapping function extends the PINN network’s capability to represent information over a larger input range due to the Fourier mapping’s frequency domain expandability. In this context, sine and cosine functions, as global Fourier basis functions, are optimal. This is expressed as follows:
\begin{equation}
\zeta(\bm{x}) = 
\left[\begin{array}{c}
\cos(\bm{\Lambda} \bm{x})\\
\sin(\bm{\Lambda} \bm{x}) 
\end{array}
\right],
\label{fourier}
\end{equation}
where $\bm{\Lambda}$ is a pre-input vector or matrix, which may be trainable or non-trainable, directly influencing the first hidden layer of the DNN and matching the number of nodes in that layer. Using this Fourier approximation, the input function $\mathcal{F}(x)$ can be represented as a combined form of sine and cosine:
\begin{equation}\label{fourierpp}
\mathcal{F}(x) =  \sum_{n=1}^{N}S\bigg{(}\sin(\bm{\Lambda}_nx),\cos(\bm{\Lambda}_nx);\bm{\tilde{\theta}}\bigg{)},
\end{equation}
where $S(x,\bm{\tilde{\theta}})$ is a DNN in the form of a fully connected network or a submodule of the DNN. Here, {${\omega_0, \omega_1, \omega_2, \cdots}$ }are the target frequency ranges within the function, with $\omega=0$ included as the base. This form resembles a Fourier expansion, with the DNN learning the Fourier basis function coefficients (except in the first layer).

\begin{figure}[htbp]
	\begin{center}
		\begin{tikzpicture}[scale=0.7]			
		\node (in) at  (-4.25,4){};
		
		\node[circle, fill=blue!50,inner sep=2.25pt] (x) at  (-1.0,4){\large$\bm{x}$};
		
		\draw[line width=0.8pt,color=blue,->] (in) --node[above]{\Large{input}} (x);
		
		\node[rectangle, rounded corners=1mm,fill=gray!50,inner sep=3.5pt] (h01) at (3.5,0){\small$\begin{matrix} \sin(W^{[1]}_{Q} \tilde{\bm{x}})\\ \cos(W^{[1]}_Q \tilde{\bm{x}})\end{matrix}$};
		\node[rectangle, rounded corners=1mm,fill=gray!50,inner sep=3.5pt] (h02) at (3.5,2){\small$\begin{matrix} \sin(W^{[1]}_{q} \tilde{\bm{x}})\\ \cos(W^{[1]}_{q} \tilde{\bm{x}})\end{matrix}$};
		\node[rectangle, rounded corners=1mm,fill=gray!0,inner sep=3.5pt] (h03) at (3.5,4){$\vdots$};
		\node[rectangle, rounded corners=1mm,fill=gray!50,inner sep=3.5pt] (h04) at (3.5,6){\small$\begin{matrix} \sin(W^{[1]}_{2} \tilde{\bm{x}})\\ \cos(W^{[1]}_{2} \tilde{\bm{x}})\end{matrix}$};
		\node[rectangle, rounded corners=1mm,fill=gray!50,inner sep=3.5pt] (h05) at (3.5,8){\small$\begin{matrix} \sin(W^{[1]}_1 \tilde{\bm{x}})\\ \cos(W^{[1]}_1 \tilde{\bm{x}})\end{matrix}$};
		
		\draw[line width=0.8pt,color=brown,->] (x) -- node[below, rotate=-50,yshift=0.5mm,color=black]{\small $\tilde{\bm{x}}=a_Q\bm{x}$}(2.15,0);
		\draw[line width=0.8pt,color=brown,->] (x) -- node[above,rotate=-25,yshift=-0.5mm,color=black]{\small$\tilde{\bm{x}}=a_q\bm{x}$}(2.15,2);
		\draw[line width=0.8pt,color=brown,->] (x) -- node[below, xshift=-0.5mm,rotate=40,color=black]{\small $\tilde{\bm{x}}=a_2\bm{x}$}(2.15,6);
		\draw[line width=0.8pt,color=brown,->] (x) -- node[above, yshift=-0.5mm,rotate=50,color=black]{\small $\tilde{\bm{x}}=a_1\bm{x}$}(2.15,8);
		
		\node[circle, fill=gray!50,inner sep=1pt] (h11) at (8,0.0){\large$\int$};
		\node[circle, fill=gray!50,inner sep=1pt] (h12) at (8,2){\large$\int$};
		\node[circle, fill=gray!0,inner sep=3pt] (h13) at (8,4){\large$\vdots$};
		\node[circle, fill=gray!50,inner sep=1pt] (h14) at (8,6){\large$\int$};
		\node[circle, fill=gray!50,inner sep=1pt] (h15) at (8,8){\large$\int$};
		
		\draw[line width=0.8pt,color=brown,->] (4.9,0) -- (h11);
		\draw[line width=0.8pt,color=brown,->] (4.9,0) -- (h12);
		\draw[line width=0.8pt,color=brown,->] (4.9,0) -- (h14);
		\draw[line width=0.8pt,color=brown,->] (4.9,0) -- (h15);
		
		\draw[line width=0.8pt,color=brown,->] (4.9,2) -- (h11);
		\draw[line width=0.8pt,color=brown,->] (4.9,2) -- (h12);
		\draw[line width=0.8pt,color=brown,->] (4.9,2) -- (h14);
		\draw[line width=0.8pt,color=brown,->] (4.9,2) -- (h15);
		
		\draw[line width=0.8pt,color=brown,->] (4.9,6) -- (h11);
		\draw[line width=0.8pt,color=brown,->] (4.9,6) -- (h12);
		\draw[line width=0.8pt,color=brown,->] (4.9,6) -- (h14);
		\draw[line width=0.8pt,color=brown,->] (4.9,6) -- (h15);
		
		\draw[line width=0.8pt,color=brown,->] (4.9,8) -- (h11);
		\draw[line width=0.8pt,color=brown,->] (4.9,8) -- (h12);
		\draw[line width=0.8pt,color=brown,->] (4.9,8) -- (h14);
		\draw[line width=0.8pt,color=brown,->] (4.9,8) -- (h15);
		
		\node[circle, fill=gray!50,inner sep=1pt] (h21) at (10.75,0){\large$\int$};
		\node[circle, fill=gray!50,inner sep=1pt] (h22) at (10.75,2){\large$\int$};
		\node[circle, fill=gray!0,inner sep=1pt] (h23) at (10.75,4){$\vdots$};
		\node[circle, fill=gray!50,inner sep=1pt] (h24) at (10.75,6){\large$\int$};
		\node[circle, fill=gray!50,inner sep=1pt] (h25) at (10.75,8){\large$\int$};
		
		\node[circle, fill=red!60,inner sep=2pt] (u) at (12.5,4){$\int$};
  
		\node (out) at (16,4){};
		
		\draw[line width=0.8pt,color=brown,->] (h11) -- (h21);
		\draw[line width=0.8pt,color=brown,->] (h11) -- (h22);
		\draw[line width=0.8pt,color=brown,->] (h11) -- (h24);
		\draw[line width=0.8pt,color=brown,->] (h11) -- (h25);
		
		\draw[line width=0.8pt,color=brown,->] (h12) -- (h21);
		\draw[line width=0.8pt,color=brown,->] (h12) -- (h22);
		\draw[line width=0.8pt,color=brown,->] (h12) -- (h24);
		\draw[line width=0.8pt,color=brown,->] (h12) -- (h25);
		
		\draw[line width=0.8pt,color=brown,->] (h14) -- (h21);
		\draw[line width=0.8pt,color=brown,->] (h14) -- (h22);
		\draw[line width=0.8pt,color=brown,->] (h14) -- (h24);
		\draw[line width=0.8pt,color=brown,->] (h14) -- (h25);
		
		\draw[line width=0.8pt,color=brown,->] (h15) -- (h21);
		\draw[line width=0.8pt,color=brown,->] (h15) -- (h22);
		\draw[line width=0.8pt,color=brown,->] (h15) -- (h24);
		\draw[line width=0.8pt,color=brown,->] (h15) -- (h25);
		
		\draw[line width=0.8pt,color=brown,->] (h21) -- (u);
		\draw[line width=0.8pt,color=brown,->] (h22) -- (u);
		\draw[line width=0.8pt,color=brown,->] (h24) -- (u);
		\draw[line width=0.8pt,color=brown,->] (h25) -- (u);
		
		\draw[line width=0.8pt,color=blue,->] (u) --node[above]{\Large{output}} (out);
		\draw[thick, dashed, draw = purple](-1.75, 9.25)-- (13.2, 9.25) -- (13.2, -1.25) -- (-1.75, -1.25) -- (-1.75, 9.25);					
		\end{tikzpicture}
	\end{center}
	\caption{\small The figure shows a Fourier basis transformation in a neural network with three hidden layers. The weight of the $q^{\text{th}}$ hidden unit in the first hidden layer is denoted by $W^{[1]}_q$.}
	\label{fourier_basis}
\end{figure}
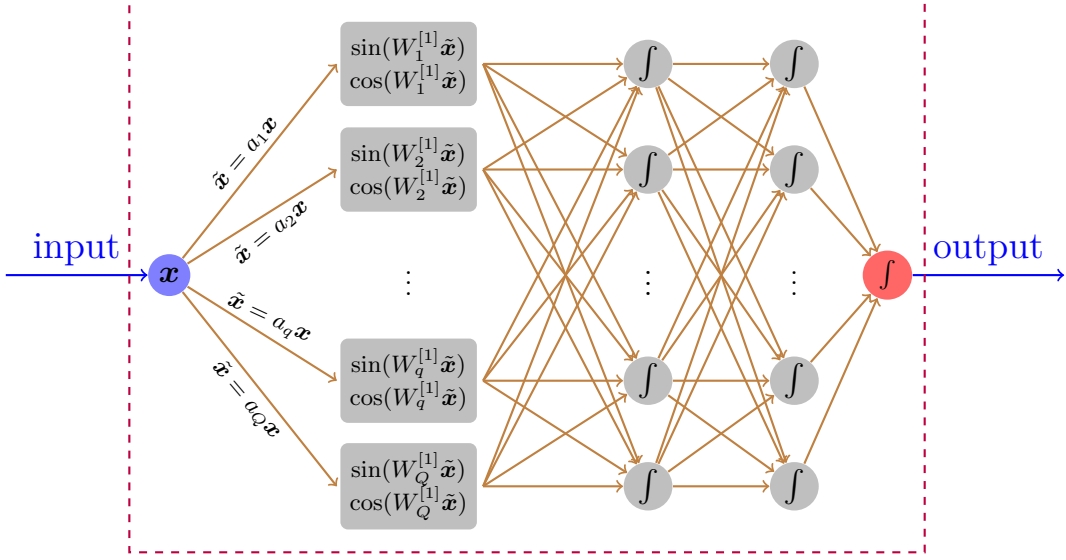

In analyzing the mechanism behind Fourier conversion and decomposition, the initial layer of the DNN model can be interpreted as a set of basis elements in Fourier space, with its output representing a linear combination of these basis functions \cite{liu2020multi,li2020elliptic,wang2020eigenvector}. This approach, known as Fourier feature mapping, enables the model to effectively capture and represent high-frequency components by leveraging the inherent properties of Fourier bases, thus enhancing performance on multi-scale and high-frequency problems.

\subsection{FCPINN Algorithm Process}
To summarize, the process of FCPINN for solving biharmonic equations is as follows:

\begin{algorithm}[H]
	\caption{FCPINN for Solving Biharmonic Equations}
	1. Construct the $k^{th}$ training set $\mathcal{S}^k$, comprising interior points $S_{I}^k=\{\bm{x}^i_I\}_{i=1}^{N_{in}}$ and boundary points $S_{B}^k=\{\bm{x}^j_B\}_{j=1}^{N_{bd}}$. Here, the point sets $\bm{x}_I^i$ and $\bm{x}_B^j$ are sampled randomly from $\mathbb{R}^d$ using the Latin hypercube sampling method.
	
	2. Compute the total objective function $\mathcal{L}(\mathcal{S}^k;\bm{\theta}^{k})$ for the given training set $\mathcal{S}^k$:
    \begin{equation*}
    \mathcal{L}(\mathcal{S}^k;\bm{\theta}^{k}) = \mathcal{L}_{in}(S_{I}^k;\bm{\theta}^k) + \gamma\mathcal{L}_{bd}(S_{B}^k;\bm{\theta}^k)
    \end{equation*}
   where $\mathcal{L}_{in}(\cdot;\bm{\theta}^k)$ is defined in \eqref{loss2ellptic}, and $\mathcal{L}_{bd}(\cdot;\bm{\theta}^k)$ is defined in \eqref{loss_bd2dirichlet} or \eqref{loss_bd2navier}.
	
	3. Update the DNN’s internal parameters using a suitable optimization method; for example, with the SGD method:
	\begin{equation*}
	\bm{\theta}^{k+1}=\bm{\theta}^{k}-\alpha^k\nabla_{\bm{\theta}^k}\mathcal{L}(\bm{x}^k;\bm{\theta}^{k})~~\text{where}~~\bm{x}^k\in\mathcal{X}^k,
	\end{equation*}
 where $\alpha^k$ (the learning rate) decreases as $k$ increases.
	
	4. Repeat steps 1-3 until the relative error meets a specified convergence criterion or the objective function stabilizes.
\end{algorithm}	

\begin{remark}
In practice, the scale factor $\bm{\Lambda}$ in \eqref{fourierpp} is a user-defined positive vector, either trainable or fixed, with a length smaller than the number of neural units in the first hidden layer. The scale factor vector is then expanded to match the size of the first hidden layer.
\end{remark}

\section{Numerical experiments}\label{sec:04}

In this section, we consider biharmonic equations of various dimensions as given in \eqref{eq:biharmonic}, with boundary constraints from either \eqref{eq:Dirichlet} or \eqref{eq:Navier}, to assess the feasibility and effectiveness of our FCPINN model. Standard PINN and MIM models with a second-order system scheme, along with their Fourier-enhanced versions, are introduced as baselines for comparison. The specifics of each model used in the numerical experiments are outlined as follows:

\begin{itemize}
    \item \textbf{PINN:} The solver is configured as a standard DNN model, with sine as the activation function for all hidden layers and a linear output layer.
    \item \textbf{MIM:} This solver is also configured as a standard DNN model, using sine as the activation function for all hidden layers and a linear output layer.
    \item \textbf{FCPINN:} The solver is configured as a standard DNN model. The Fourier feature mapping is applied as the activation function in the first hidden layer, while sine functions are used in the remaining hidden layers. The output layer is linear, with the scale factor chosen as $\Lambda=(1,2,3,...,14,15)$.
\end{itemize}

Additionally, we examine the impact of different activation functions on our proposed model, including ReLU, ReQU, GELU, sigmoid, tanh, and sine. For value ranges, the problem domain and boundary are located in Euclidean space $\mathbb{R}^d$. Training and testing data are sampled uniformly from the target region. All neural networks are trained using the Adam optimizer, with an initial learning rate of 0.01 that decays by 2.5\% every 100 iterations \cite{kingma2015adam}. Model performance is evaluated using the following criterion:
\begin{equation*}
REL = \sqrt{\sum_{i=1}^{N}\frac{|\tilde{u}(\bm{x}^i)-u^*(\bm{x}^i)|^2}{|u^*(\bm{x}^i)|^2}},
\end{equation*}
where $\tilde{u}(\bm{x}^i)$ and $u^*(\bm{x}^i)$ represent the DNN-predicted and true values at testing points $\{\bm{x}^i\}(i=1,2,\cdots,N)$, respectively. Here, $N$ denotes the number of test sample points. During training, sampling and testing occur every 1000 iterations to monitor the model’s performance in real-time. The penalty coefficient $\beta$ is set according to the following rule:
\begin{equation}
\beta=\left\{
\begin{aligned}
\beta_0, \quad &\textup{if}~~i<M_{\max}*0.1\\
10\beta_0,\quad &\textup{if}~~M_{\max}*0.1 \leq i < M_{\max}*0.2\\
50\beta_0, \quad&\textup{if}~~ M_{\max}*0.2 \leq i < M_{\max}*0.25\\
100\beta_0, \quad&\textup{if}~~ M_{\max}*0.25 \leq i < M_{\max}*0.5\\
200\beta_0, \quad&\textup{if}~~ M_{\max}*0.5 \leq i < M_{\max}*0.75\\
500\beta_0, \quad&\textup{otherwise},
\end{aligned}
\right.
\end{equation}
where $\beta_0=10$ for all tests, and $M_{\max}$ represents the maximum number of epochs.

All training and testing processes are conducted in PyTorch (version 1.12.1). The experimental platform is a workstation with 64 GB RAM and a single NVIDIA GeForce RTX 4090 24-GB GPU.

\subsection{Case for Dirichlet Boundary}

\begin{example}\label{2D_Dirichlet_E1}
We first consider the biharmonic equation with Dirichlet boundary conditions in a two-dimensional Euclidean space, where the force term depends only on $x_1$ and $x_2$, indicating a linear problem. The domain of interest is a rectangular region $\Omega=[x_1^{min}, x_1^{max}]\times[x_2^{min}, x_2^{max}]$, and an exact solution is given by
\begin{equation*}
u(x_1,x_2)=e^{(x_1-x_1^{min})(x_1^{max}-x_1)(x_2-x_2^{min})(x_2^{max}-x_2)}.
\end{equation*}
This naturally induces Dirichlet boundary functions $g(x_1,x_2)$ and $h(x_1,x_2)$, defined as
\begin{equation*}
    \begin{cases}
        g(x_{1}^{min}, x_2) = g(x_{1}^{max}, x_2) = g(x_{1}, x_{2}^{min}) = g(x_{1}, x_{2}^{max}) = 1,\\
        h(x_1^{min},x_2)=-(x_1^{min}+x_1^{max}-2x_1)(x_2-x_2^{min})(x_2^{max}-x_2),\\
        h(x_1^{max},x_2)=(x_1^{min}+x_1^{max}-2x_1)(x_2-x_2^{min})(x_2^{max}-x_2),\\
        h(x_1,x_2^{min})=-(x_1-x_1^{min})(x_1^{max}-x_1)(x_2^{min}+x_2^{max}-2x_2),\\
        h(x_1,x_2^{max})=(x_1-x_1^{min})(x_1^{max}-x_1)(x_2^{min}+x_2^{max}-2x_2).\\
    \end{cases} 
\end{equation*}
With careful calculations, the force term $f(x_1, x_2)$ can be obtained, though it is omitted here.
\end{example}

Initially, we approximate the solution of \eqref{eq:biharmonic} using the three methods described on the regular domain $[-1,1]\times[-1,1]$, with each solver consisting of four hidden layers containing 30 units each. In this case, the model underwent 50,000 iterations. In each iteration, all three models were trained using 3000 sample points within the domain and 2000 sample points on the boundaries. The final training times and relative errors for each model are provided in Table 1, while Figures 5 and 6 illustrate the absolute errors of the output $u$ and the intermediate function $v$. To evaluate the training speed and accuracy, we uniformly sampled 16,384 equidistant grid points in the domain, and the real-time relative error (REL) of the three models during training is shown in Figure 5(e), with the error trend recorded every 1000 iterations.

First, from a general perspective, as shown in Table 1, the accuracy of FCPINN and FPINN is comparable and significantly higher than that of MIM. In terms of computational speed, FMIM and FCPINN perform similarly and are notably faster than FPINN. This is because FCPINN combines a coupled structure, benefiting from both MIM and FPINN, achieving an optimal balance between speed and accuracy.

Looking at the error distribution in detail, the heatmaps in Figures 5(b), 5(c), and 5(d) reveal that FCPINN achieves higher accuracy than FPINN and FMIM. This demonstrates that FCPINN, with its Fourier expansion and coupled framework, simplifies computations and effectively mitigates challenges such as high-order oscillations in DNN parameters, as well as gradient vanishing or explosion in high-order problems.

\begin{figure}[!ht]
    \centering
    \subfigure[Exact solution]{
        \label{UExact}
        \includegraphics[scale=0.34]{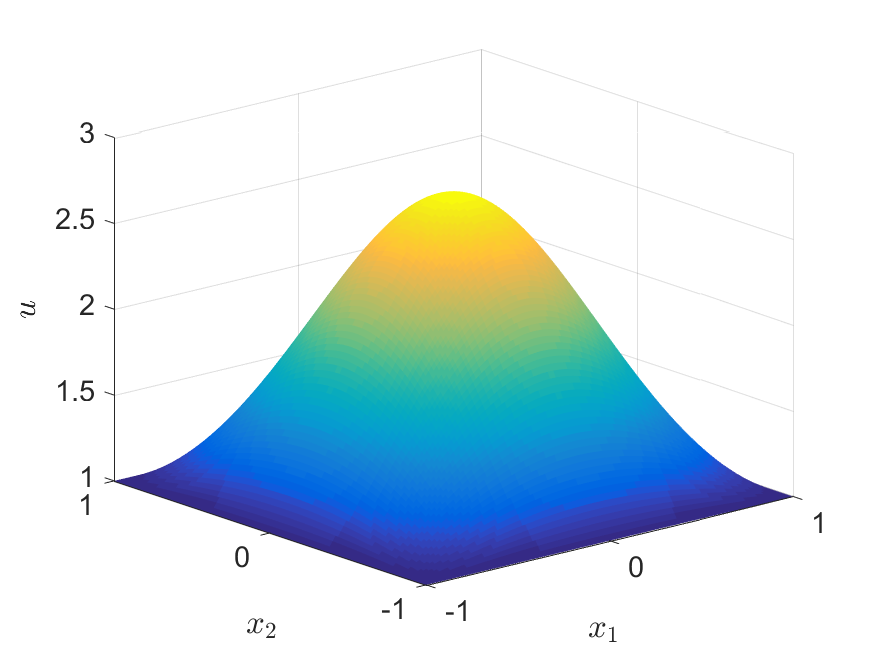}
    }
    \subfigure[Pointwise error for PINN]{
        \label{UABS2PINN_Birichlet2D_E1}
        \includegraphics[scale=0.34]{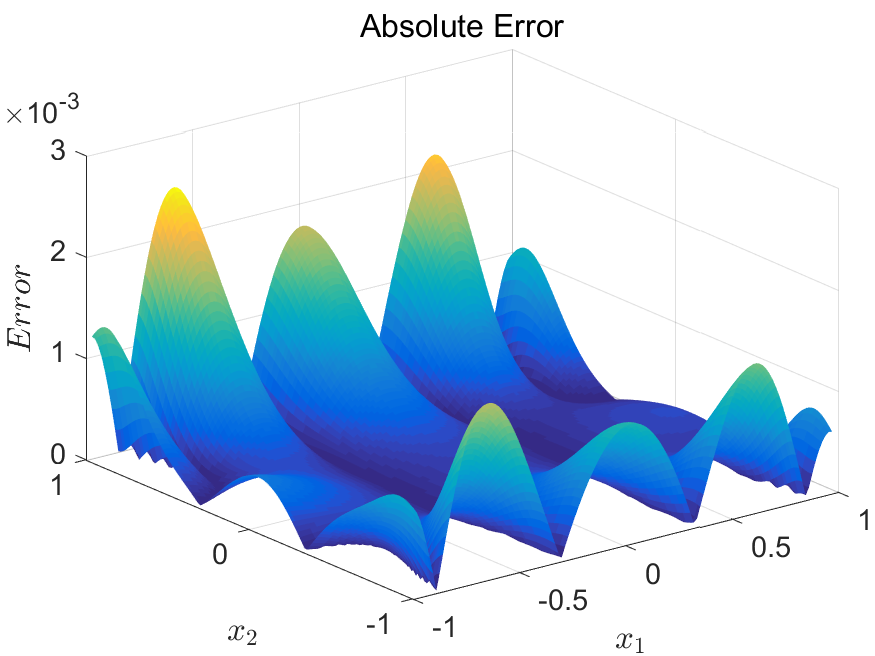}
    }
    \subfigure[Pointwise error for MIM]{
        \label{UABS2MIM_Birichlet2D_E1}
        \includegraphics[scale=0.34]{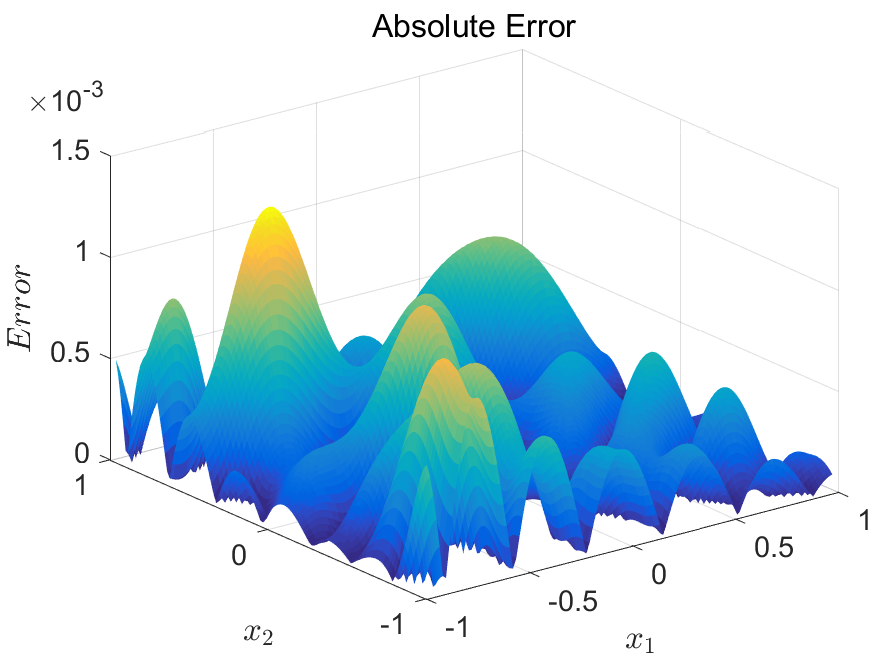}
    }
    \subfigure[Pointwise error for FCPINN]{
        \label{UABS2CPINN_Birichlet2D_E1}
        \includegraphics[scale=0.35]{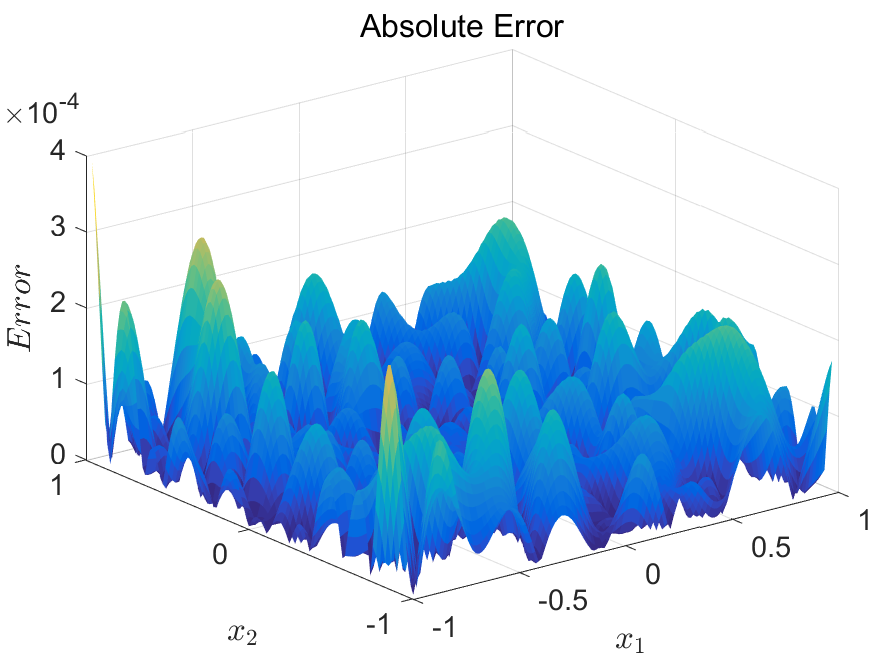}
    }
    \subfigure[REL]{
       \label{REL2Birichlet2D_E1}
    \includegraphics[scale=0.33]{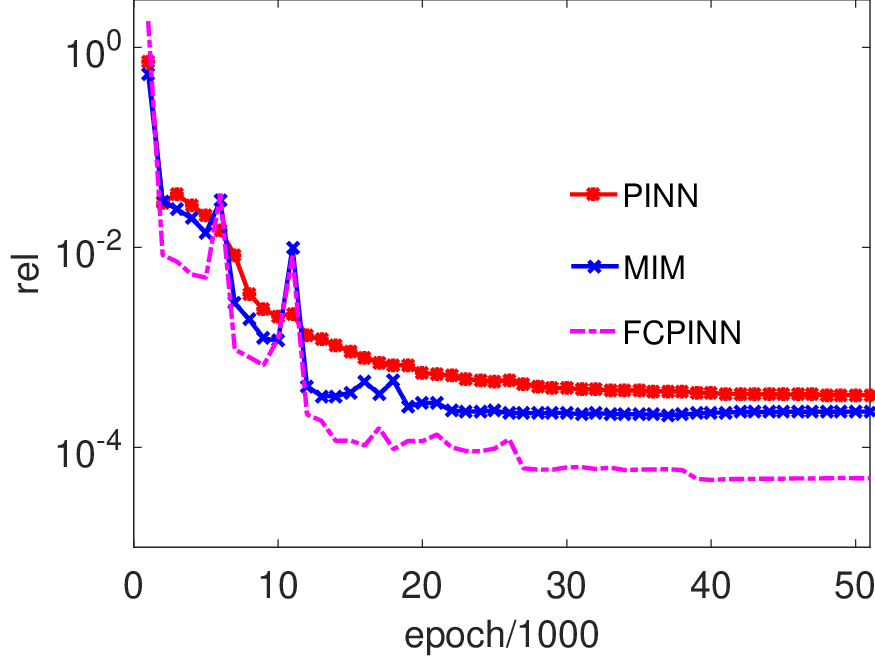}
    }
    \label{BiDirichlet_2DE1_U}
    \caption{The test results of solution for Example \ref{2D_Dirichlet_E1}.}
\end{figure}

Additionally, as indicated by the curve in Figure~\ref{REL2Birichlet2D_E1}, FCPINN shows a lower initial error and faster convergence rate than standard FPINN and FMIM, illustrating that the coupling algorithm inspired by FMIM significantly improves performance in FPINN by avoiding high-order derivatives.

In conclusion, this example shows that for high-order problems with Dirichlet boundary conditions and wide frequency ranges, FCPINN outperforms traditional FPINN and FMIM in both accuracy and training speed.

\begin{table}[!ht]
	\centering
	\caption{The REL and Time of the above three DNN models, for example, \ref{2D_Dirichlet_E1} on various domain}
	\label{Table_2D_Dirichlet_E1}
	\begin{tabular}{|l|c|c|c|}
		\hline  
		             &PINN                 & MIM                  &  FCPINN               \\ \hline
		  REL       &$3.342\times 10^{-4}$&$2.276\times 10^{-4}$ &$4.905\times 10^{-5}$  \\ 
		 Time(s)     &1938.947             &792.726              & 981.735              \\  \hline
	\end{tabular}
\end{table}

\begin{figure}[H]
    \centering
    \subfigure[Exact solution]{
        \label{VExact_E1}
        \includegraphics[scale=0.33]{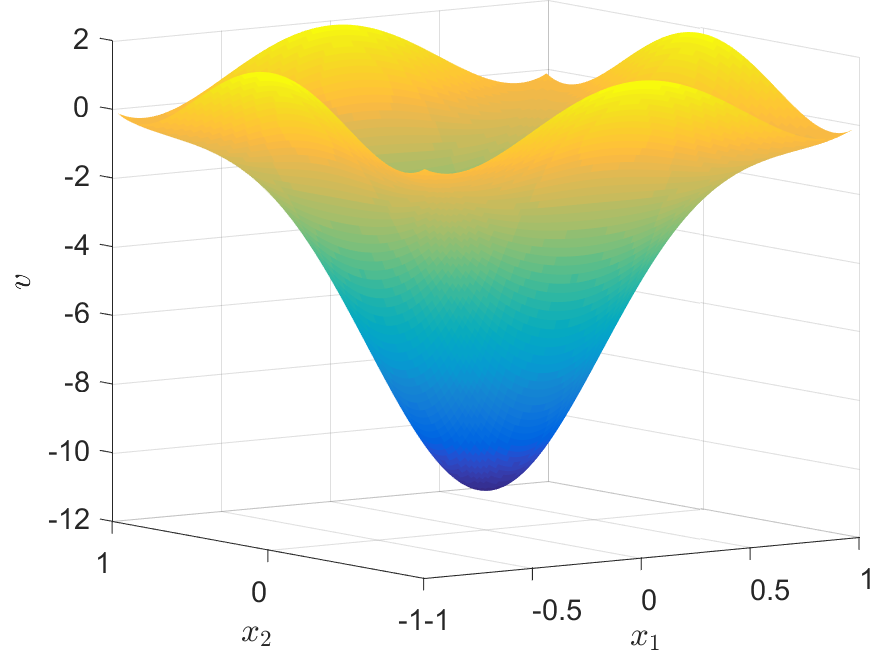}
    }
    \subfigure[Pointwise error for MIM]{
        \label{VABS2MIM_Birichlet2D_E1}
        \includegraphics[scale=0.33]{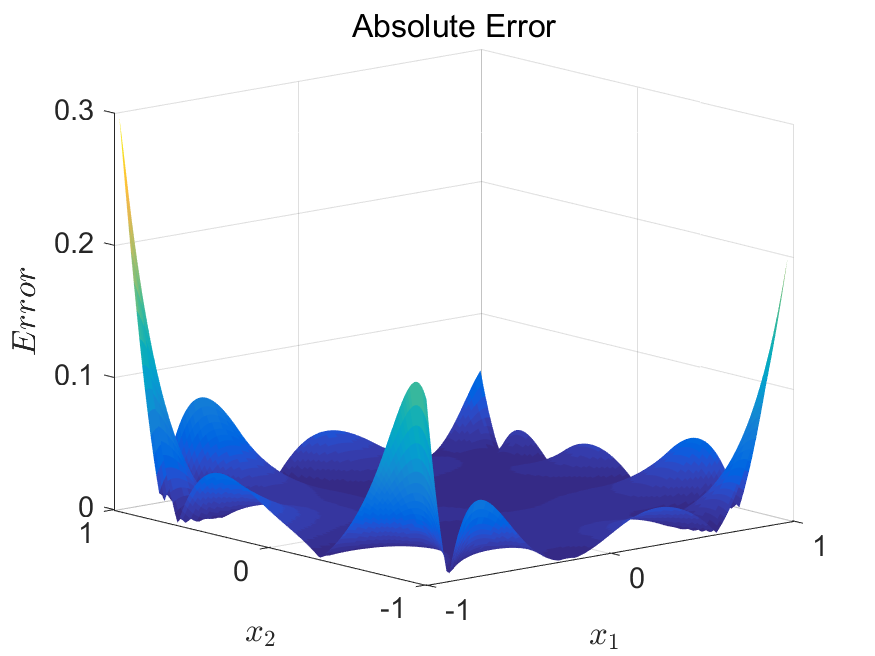}
    }
    \subfigure[Pointwise error for FCPINN]{
        \label{VABS2CPINN_Birichlet2D_E1}
        \includegraphics[scale=0.33]{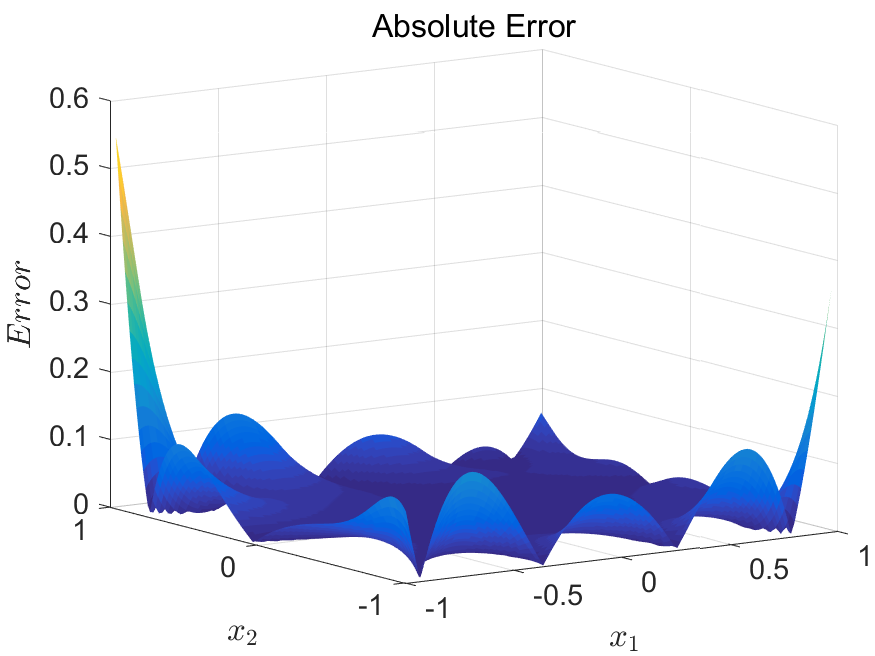}
    }
    \label{BiDirichlet_2DE1_V}
    \caption{The test results of $v$ for Example \ref{2D_Dirichlet_E1}.}
\end{figure}

\emph{Impact of Activation Function on FCPINN:} Next, we examine the effect of activation function choice on the FCPINN model for solving \eqref{eq:biharmonic}. We test sin, tanh, sigmoid, gelu, relu, and requ as activation functions for the hidden layers of each subnetwork (excluding the first and last layers), applying these functions to the PINN, MIM, and FCPINN models to solve the biharmonic problem, with all other settings as in previous experiments. Overall, we observe that the FCPINN model consistently shows lower error than the other two models, with sin and gelu demonstrating particularly competitive performance. Consequently, for subsequent models, we select $\sin$ as the activation function for our FCPINN model, and further details will not be repeated.

\begin{table}[H]
	\centering
	\caption{REL under diverse choices of activation functions with PINN, MIM and FCPINN employed for the Example \ref{2D_Dirichlet_E1}.}
	\label{Table2activation}
	\begin{tabular}{|l|c|c|c|c|c|c|c|}
	\hline
            &sin                  &tanh                &sigmoid      & gelu                & relu  & requ      \\  \hline
        PINN&$3.342\times 10^{-4}$&$2.144\times10^{-4}$&$0.0021$&$5.385\times10^{-5}$&$0.138$&$0.0075$           \\ \hline
        MIM&$2.276\times 10^{-4}$&$1.892\times10^{-4}$&$9.133\times10^{-4}$&$7.483\times10^{-5}$&$0.0949$&$0.0034$           \\ \hline
	FCPINN&$4.905\times 10^{-5}$&$4.798\times 10^{-5}$&$1.174\times10^{-4}$&$3.162\times 10^{-5}$&0.0077 &$4.031\times10^{-4}$    \\	\hline
	\end{tabular}
\end{table}

\emph{Impact of Hidden Units:} We investigated how changes in the number of nodes and layers within the hidden layers affect FCPINN model performance. Hidden layer configurations tested included $(20,20,20)$, $(30,30,30)$, $(30,30,30,30)$, and $(30,30,30,30,30)$, with all other settings consistent with previous experiments, to assess the performance differences across various hidden layer setups for the same problem. As shown in Table \ref{Table2hiddens}, increasing the number of hidden layers and nodes yields a slight improvement in FCPINN performance. However, given the exponential growth in computational resource requirements with additional layers and nodes, we selected $(30,30,30,30)$ as the configuration for our subsequent experiments, and further details will not be repeated.

\begin{table}[H]
	\centering
	\caption{REL for varying depths and widths of hidden units with FCPINN employed for the Example. \ref{2D_Dirichlet_E1}.}
	\label{Table2hiddens}
	\begin{tabular}{|l|c|c|c|}
		\hline
	    $(20,20,20)$         &$(30,30,30)$        &$(30,30,30,30)$    &$(30,30,30,30, 30)$\\ \hline
            $3.705\times 10^{-4}$&$5.744\times10^{-5}$&$4.905\times 10^{-5}$&$4.795\times 10^{-5}$    \\
		\hline
	\end{tabular}
\end{table}

\begin{example}\label{2D_Dirichlet_E2}
We consider a nonlinear case for the biharmonic equation~\eqref{eq:biharmonic}, where $f(\bm{x}, u, \Delta u) = \tilde{f}(\bm{x}) - \Delta u - u$, and solve it numerically in a 2D space within a hexagonal region $\Omega$ on the interval $[0,1]\times[0,1]$. The analytical solution is given by
\begin{equation*}
   u(x_1,x_2)=\sin(\pi x_1)\sin(\pi x_2),
\end{equation*}
which leads to
\begin{equation*}
f(x_1,x_2; u, \Delta u)=4\pi^4\sin(\pi x_1)\sin(\pi x_2) +2\pi^2\sin(\pi x_1)\sin(\pi x_2) - \sin(\pi x_1)\sin(\pi x_2),
\end{equation*}
with $\tilde{f}(\bm{x}) = 4\pi^4\sin(\pi x_1)\sin(\pi x_2)$. The boundary constraint functions $g(x_1,x_2)$ and $h(x_1,x_2)$ on $\partial\Omega$ can be easily obtained by direct computation, but are omitted here. The network configuration for this problem is identical to that used in the linear case.
\end{example}

The experimental results below indicate that, for complex nonlinear problems, the behavior of the three models is similar to that observed in the previous example. Here, the relationships between functions are implicitly represented, and due to the function complexity, PINN exhibits significantly slower computation speed and lower accuracy compared to the other two models. While the MIM method outperforms PINN in both speed and accuracy, it is still less competitive than FCPINN. This is because, for complex problems, the computational load of PINN is too high, leading to substantial degradation in accuracy and speed after multiple derivatives. Both MIM and FCPINN utilize decomposition methods to simplify calculations, but FCPINN holds a distinct advantage in handling nonlinear problems over MIM. Overall, FCPINN demonstrates superior performance to both PINN and MIM for two-dimensional problems with implicit relationship constraints.

\begin{table}[!ht]
	\centering
	\caption{The REL and Time of the above three models, for example, \ref{2D_Dirichlet_E2} on various domain}
	\label{Table_2D_Dirichlet_E2}
	\begin{tabular}{|l|c|c|c|}
		\hline  
		             & PINN                 & MIM                   &  FCPINN                 \\ \hline
		  REL       &$0.0039$&$9.028\times 10^{-4}$ &$2.236\times 10^{-4}$  \\ 
		 Time(s)     &2214.710      &1529.943               &  1828.699                     \\  \hline
	\end{tabular}
\end{table}

\begin{figure}[H]
    \centering
    \subfigure[Exact solution]{
        \label{UExact2Dirichlet_E2}
        \includegraphics[scale=0.34]{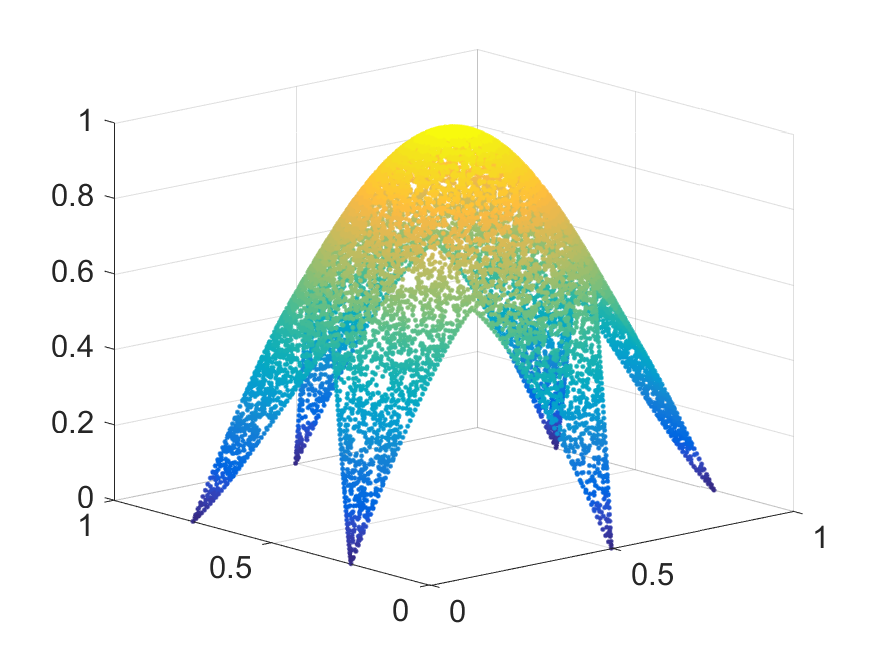}
    }
    \subfigure[Pointwise error for PINN]{
        \label{UABS2PINN_Birichlet2D_E2}
        \includegraphics[scale=0.34]{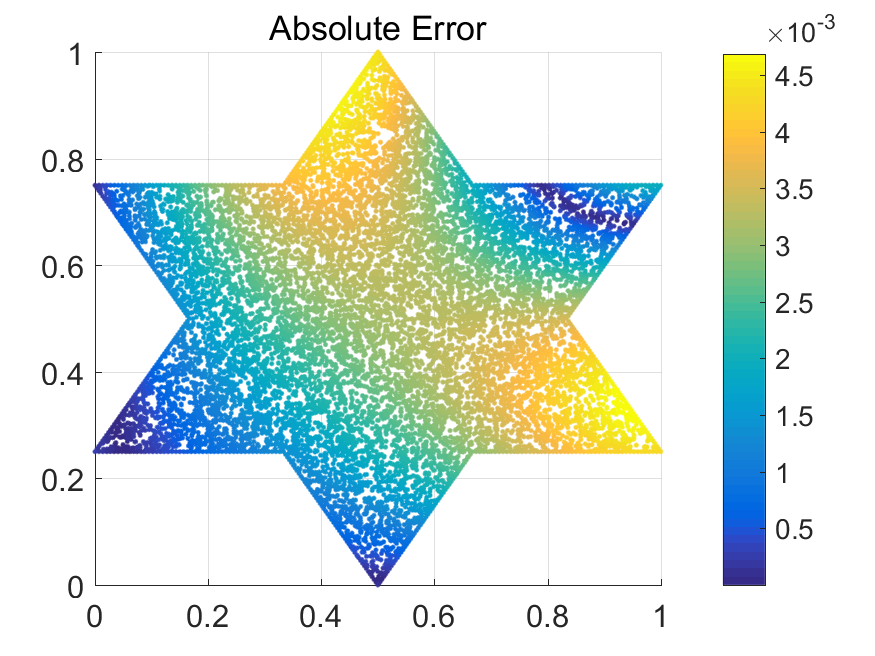}
    }
    \subfigure[Pointwise error for MIM]{
        \label{UABS2MIM_Birichlet2D_E2}
        \includegraphics[scale=0.34]{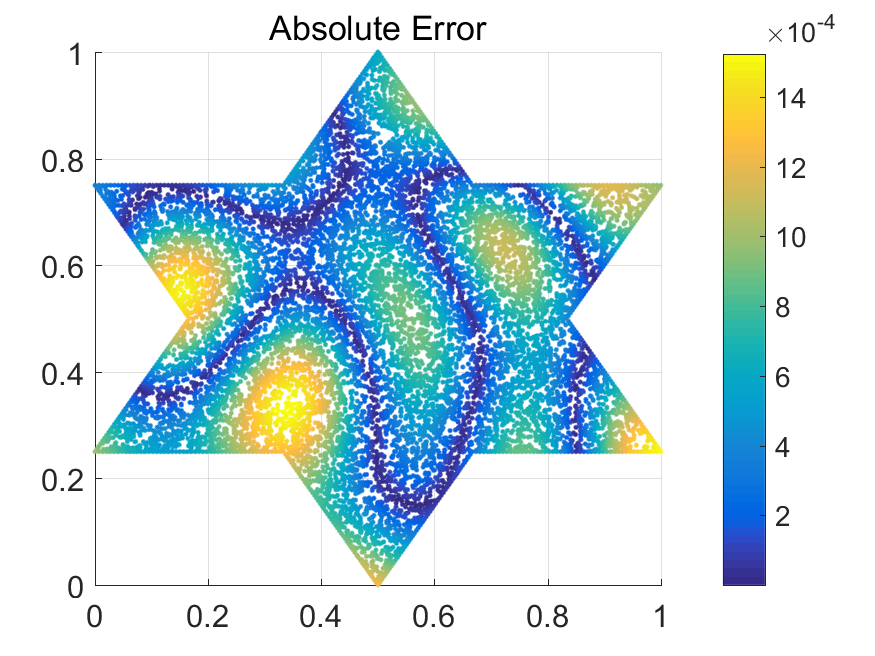}
    }
    \subfigure[Pointwise error for FCPINN]{
        \label{UABS2CPINN_Birichlet2D_E2}
        \includegraphics[scale=0.34]{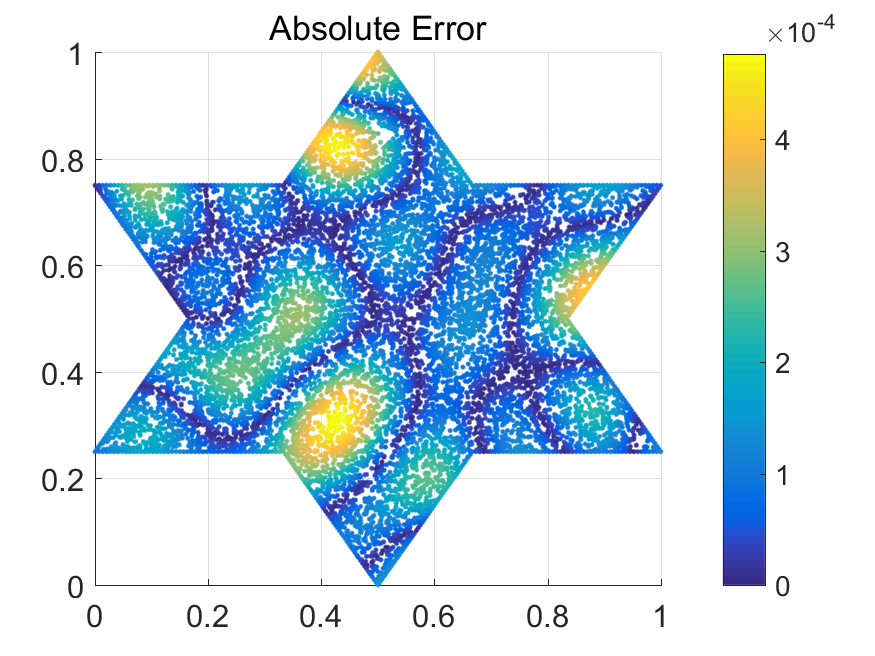}
    }
    \subfigure[REL]{
        \label{REL2Birichlet2D_E2}
        \includegraphics[scale=0.34]{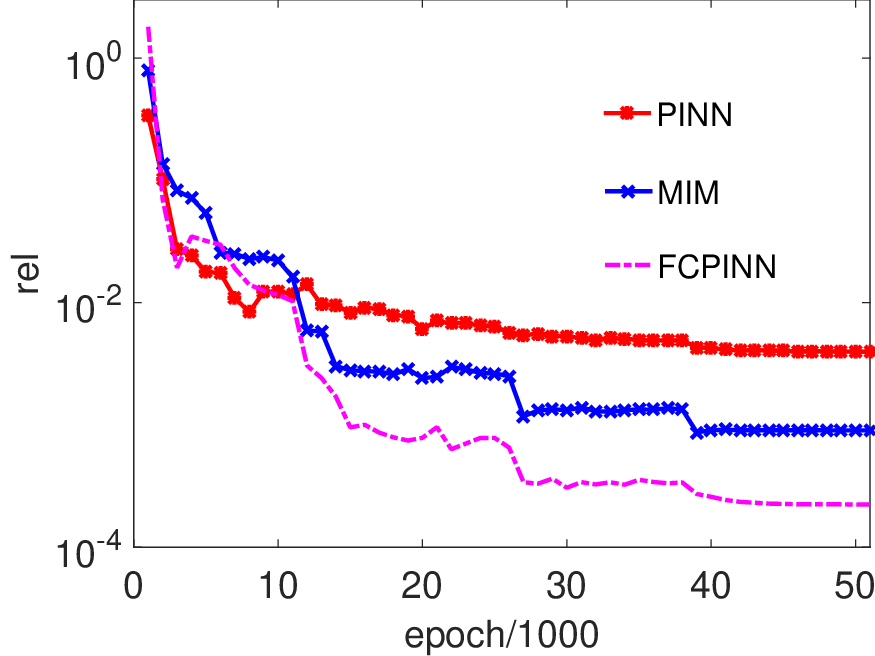}
    }
    \label{BiDirichlet_2DE2_U}
    \caption{The test results of solution for Example \ref{2D_Dirichlet_E2}.}
\end{figure}

\begin{figure}[H]
    \centering
    \subfigure[Exact solution]{
        \label{VExact_E2}
        \includegraphics[scale=0.325]{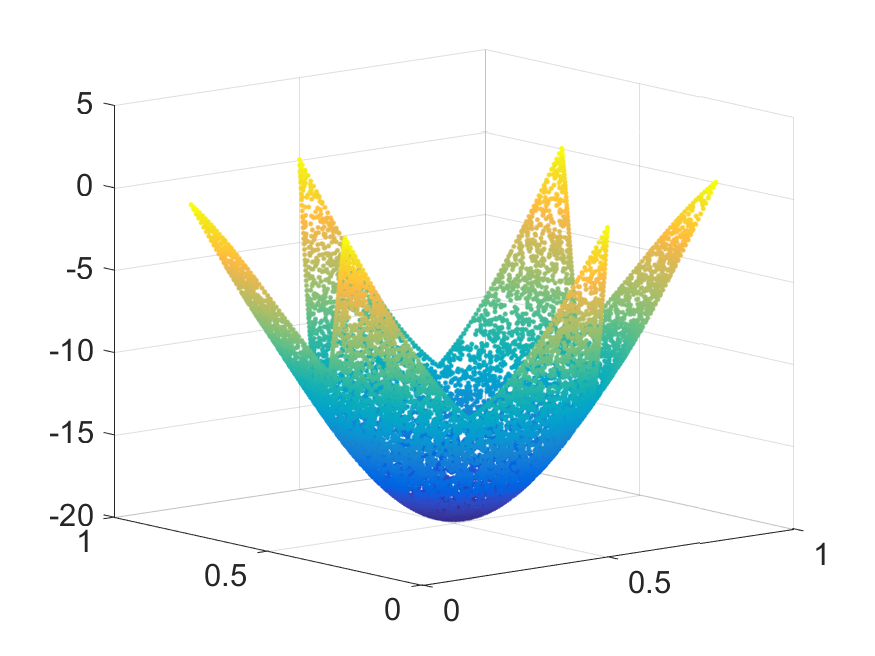}
    }
    \subfigure[Pointwise error for MIM]{
        \label{VABS2MIM_Birichlet2D_E2}
        \includegraphics[scale=0.325]{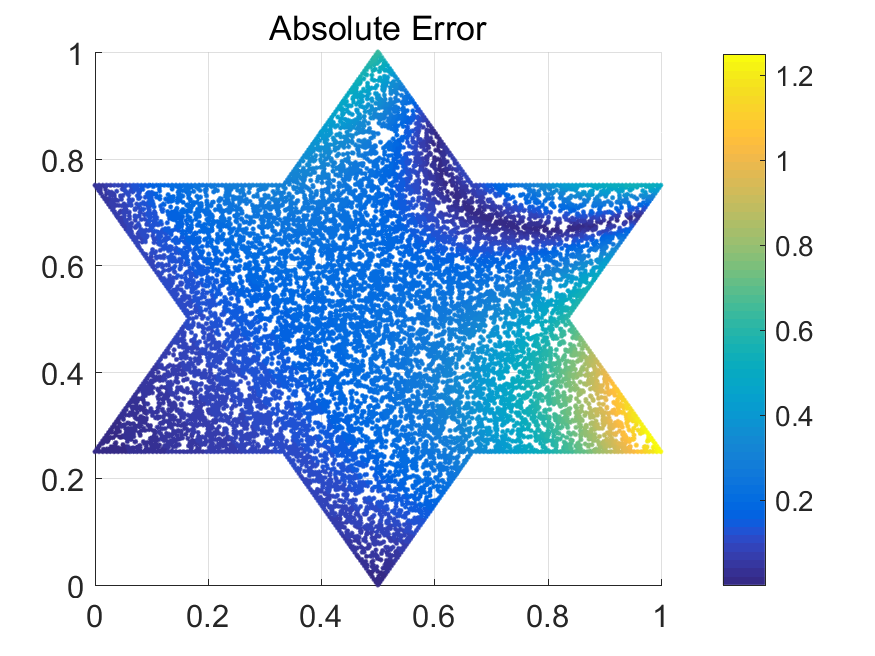}
    }
    \subfigure[Pointwise error for FCPINN]{
        \label{VABS2CPINN_Birichlet2D_E2}
        \includegraphics[scale=0.325]{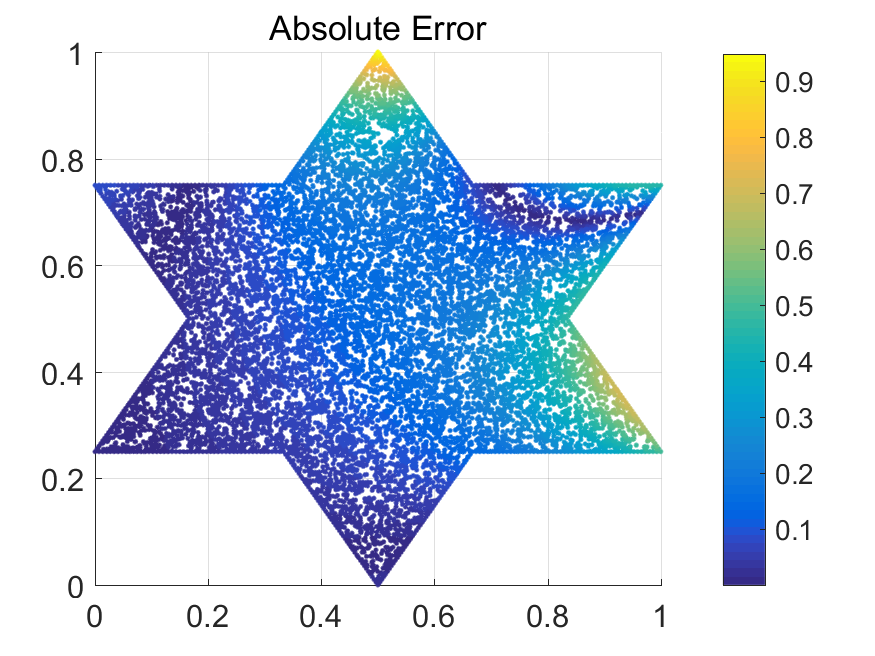}
    }
    \label{BiDirichlet_2DE2_V}
    \caption{The test results of $v$ for Example \ref{2D_Dirichlet_E2}.}
\end{figure}

\begin{example}\label{3D_Dirichlet_E1}
We now approximate the solution of a \textbf{nonlinear} biharmonic equation~\eqref{eq:biharmonic} on a complex surface within the three-dimensional cubic domain $\Omega=[-1,1]\times [-1,1]\times[-1,1]$. The analytical solution and corresponding force term are given by
\begin{equation*}
u(x_1,x_2,x_3)=50e^{-0.25(x_1+x_2+x_3)},
\end{equation*}
and
\begin{equation*}
f(x_1,x_2, x_3)=\frac{225}{128}e^{-0.25(x_1+x_2+x_3)},
\end{equation*}
respectively. The boundary functions $g(x_1,x_2,x_3)$ and $h(x_1,x_2,x_3)$ can be derived from the exact solution but are omitted here for brevity. The network configuration for this problem is the same as that used in the linear case.
\end{example}

For the 3D problem, the error ranking across the three models is \(PINN > MIM > FCPINN\). Given the increased computational load, PINN’s computation time is significantly higher than that of the other two models. While FCPINN and MIM are similar in performance, MIM is more time-efficient, whereas FCPINN demonstrates better accuracy. Notably, MIM also achieves a faster convergence rate during the initial iterations compared to both PINN methods. In conclusion, FCPINN outperforms both PINN and MIM in 3D problems with complex boundary conditions.

\begin{table}[!ht]
    \centering
    \caption{The REL and Time of the above four PIELM models, for example, \ref{3D_Dirichlet_E1} on various domain}
    \label{Table_3D_Dirichlet}
    \begin{tabular}{|l|c|c|c|}
        \hline  
                  & PINN     & MIM                  &  FCPINN               \\ \hline
          REL     &$0.1833$  &$8.612\times 10^{-4}$ &$2.828\times 10^{-5}$  \\ 
         Time(s)  & 3149.269  &1028.204              & 1223.232              \\  \hline
    \end{tabular}
\end{table}

\begin{figure}[H]
    \centering
    \subfigure[Exact solution]{
        \label{UExact2Dirichlet_3D}
        \includegraphics[scale=0.345]{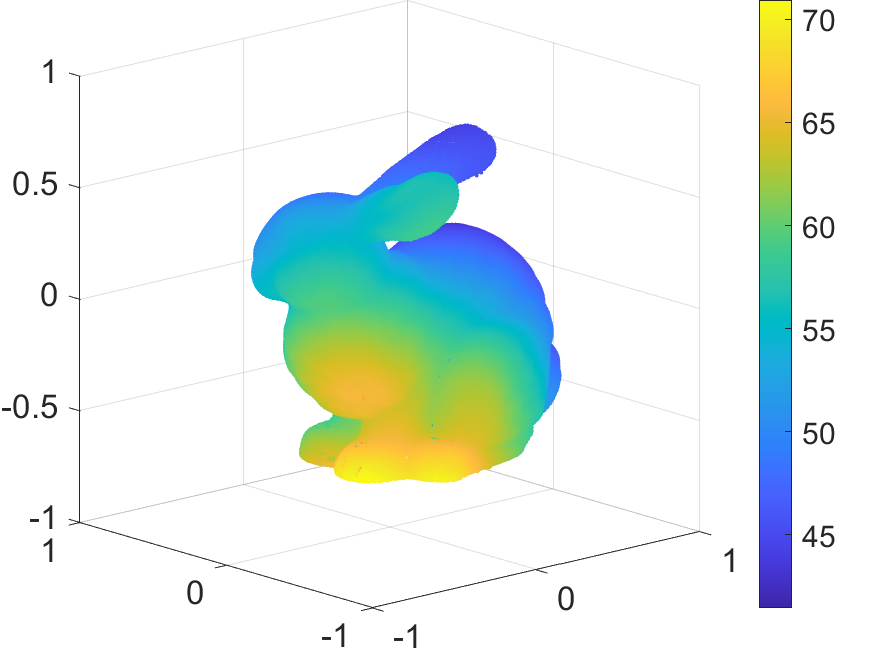}
    }
    \subfigure[Pointwise error for PINN]{
        \label{UABS2PINN_Dirichlet3D}
        \includegraphics[scale=0.345]{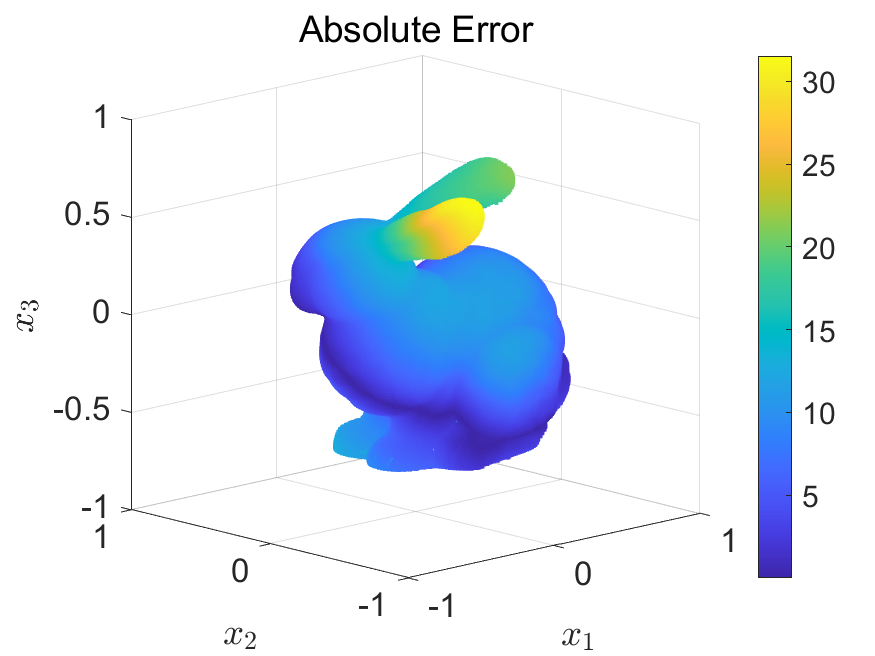}
    }
    \subfigure[Pointwise error for MIM]{
        \label{UABS2MIM_Dirichlet3D}
        \includegraphics[scale=0.345]{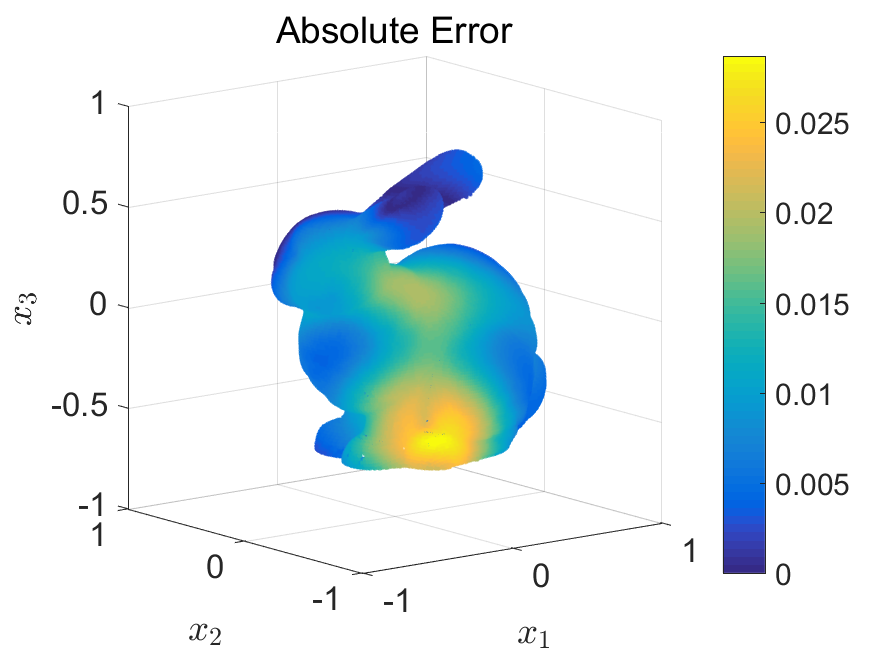}
    }
    \subfigure[Pointwise error for FCPINN]{
        \label{UABS2CPINN_Dirichlet3D}
        \includegraphics[scale=0.345]{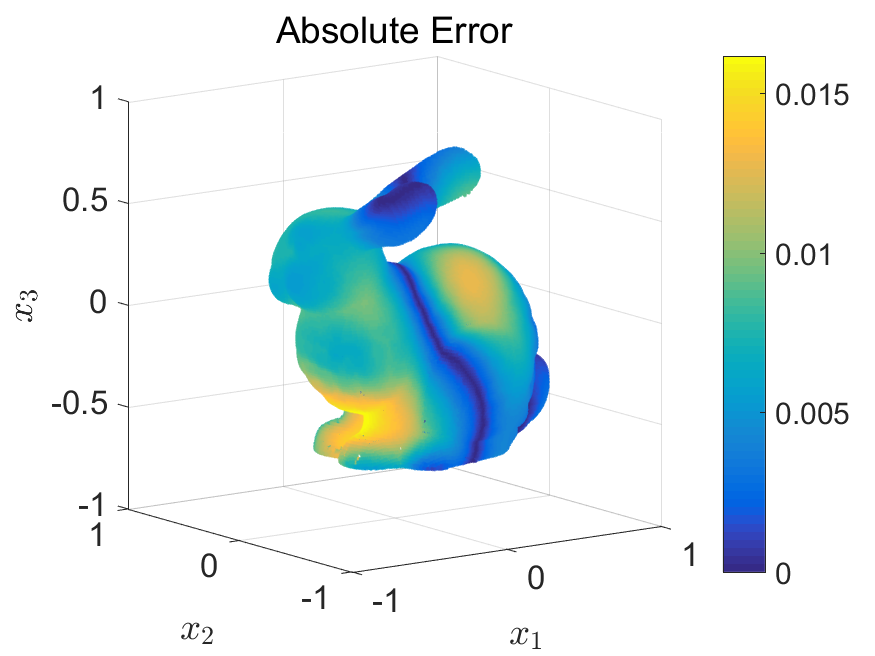}
    }
        \subfigure[REL]{
        \label{REL2Birichlet3D}
        \includegraphics[scale=0.34]{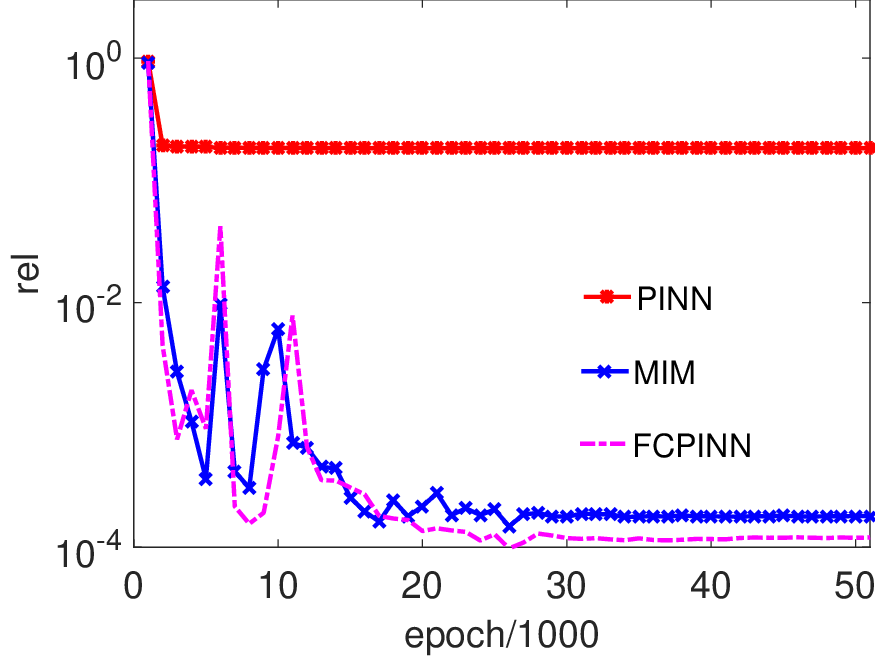}
    }
    \label{BiDirichlet_3D}
\end{figure}

\begin{figure}[H]
	\centering
	\subfigure[Exact solution]{
		\label{VExact2Dirichlet3D}
		\includegraphics[scale=0.325]{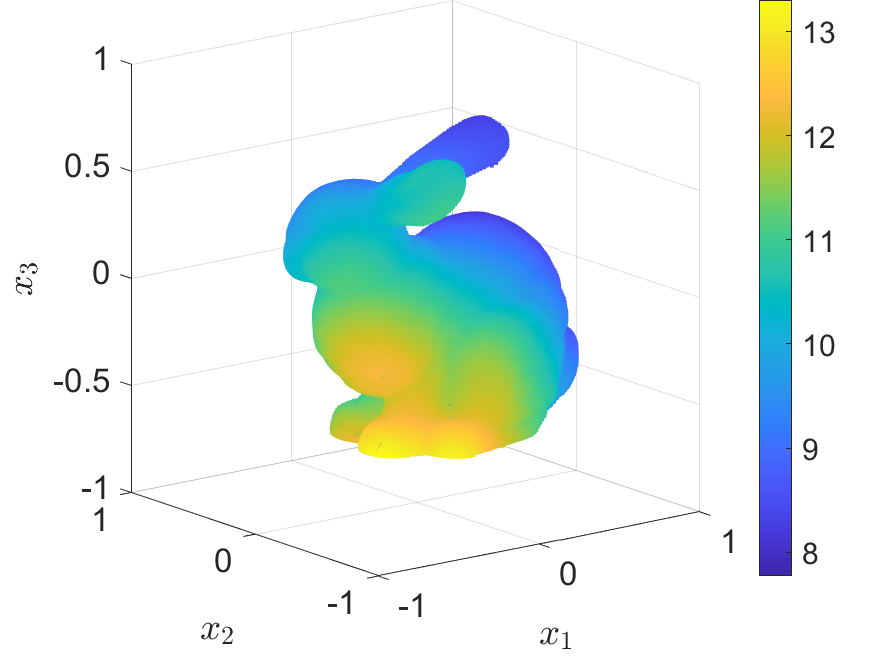}
	}
	\subfigure[Pointwise error for MIM]{
		\label{VABS2MIM_Birichlet3D}
		\includegraphics[scale=0.325]{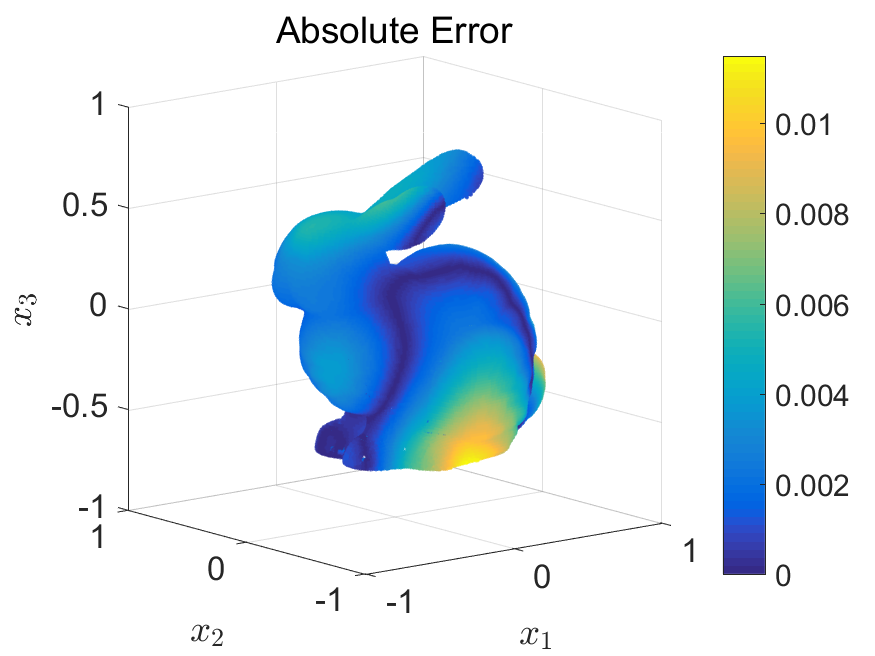}
	}
	\subfigure[Pointwise error for FCPINN]{
		\label{VABS2CPINN_Birichlet3D}
		\includegraphics[scale=0.325]{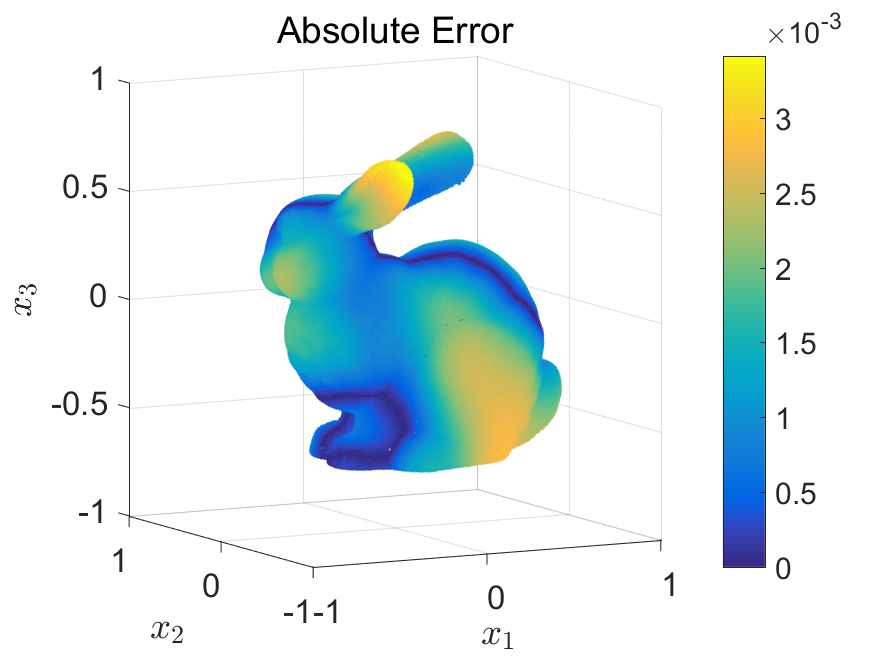}
	}
	\label{BiDirichlet3D_V}
	\caption{The test results of $v$ for Example \ref{3D_Dirichlet_E1}.}
\end{figure}

\subsection{Case for Navier Boundary}
We now examine the capability of our model to solve the biharmonic equation with Navier boundary conditions in Euclidean spaces $\mathbb{R}^2$, $\mathbb{R}^3$, and $\mathbb{R}^8$.

\begin{example}\label{BiNavier2D_E01}
We solve the biharmonic equation \eqref{eq:biharmonic} on a rectangular domain $\Omega=[0,1]\times[0,\pi]$. The exact solution and force term are given by
\begin{equation*}
u(x_1,x_2)=e^{x_1}\sin(x_2)
\end{equation*}
such that $f(x_1,x_2) = 0$ and $k(x_1,x_2)=0$. The boundary function $g(x_1, x_2)$ on $\partial \Omega$ can be derived directly from the exact solution, but we omit it here. The network configuration for this problem is the same as in the Dirichlet cases above.
\end{example}

For Navier boundary problems, the performance ranking of the three models differs from that on Dirichlet boundaries, with error ranking as \(MIM > PINN > FCPINN\), as shown in the table below. In terms of computation time, the models retain the same trend: MIM is the fastest, FCPINN is close behind, and PINN takes considerably longer. Unlike Dirichlet boundaries, Navier boundaries impose constraints on a linear combination of the function’s value and its derivative, rather than a simple boundary condition. This weaker enforcement of boundary conditions in MIM calculations leads to a significant increase in error after multiple residual evaluations. Both types of PINNs adhere more strictly to the physical boundary constraints, which results in slower computations but better preservation of the function's characteristics and thus higher accuracy.

\begin{table}[!ht]
    \centering
    \caption{The REL and Time of the above four PIELM models, for example, \ref{BiNavier2D_E01} on various domain}
    \label{Table2BiNavier2D_E1}
    \begin{tabular}{|l|c|c|c|}
        \hline  
                  & PINN     & MIM                  &  FCPINN               \\ \hline
          REL     &$4.935\times10^{-4}$  &$2.844\times 10^{-4}$ &$2.769\times 10^{-5}$  \\ 
         Time(s)  &1720.186  &527.456               &676.705               \\  \hline
    \end{tabular}
\end{table}

\begin{figure}[H]
	\centering
	\subfigure[Exact solution]{
		\label{UExact2BiNavier2D_E1}
		\includegraphics[scale=0.34]{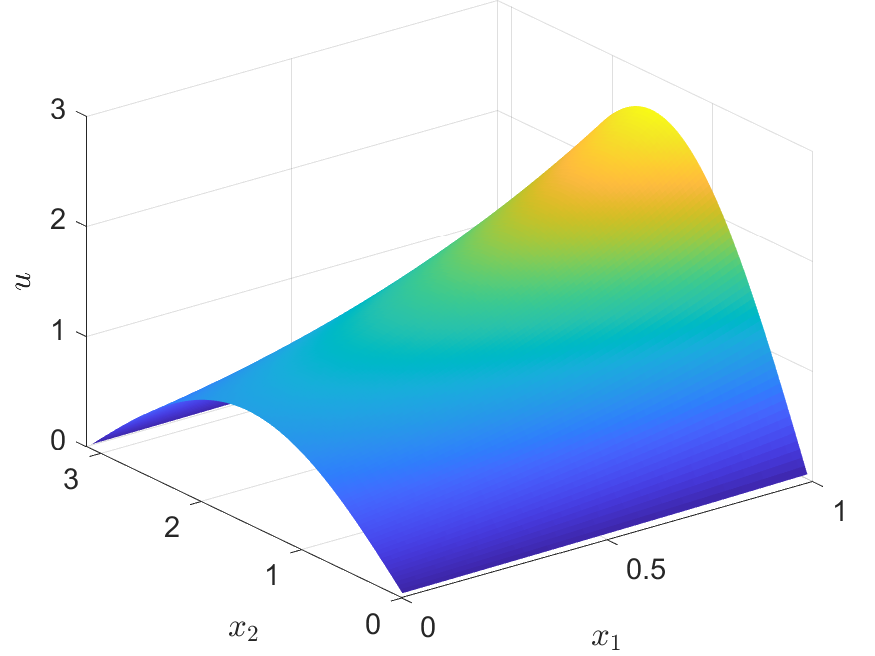}
	}
	\subfigure[Pointwise error for PINN]{
		\label{UABS2PINN_BiNavier2D_E1}
		\includegraphics[scale=0.34]{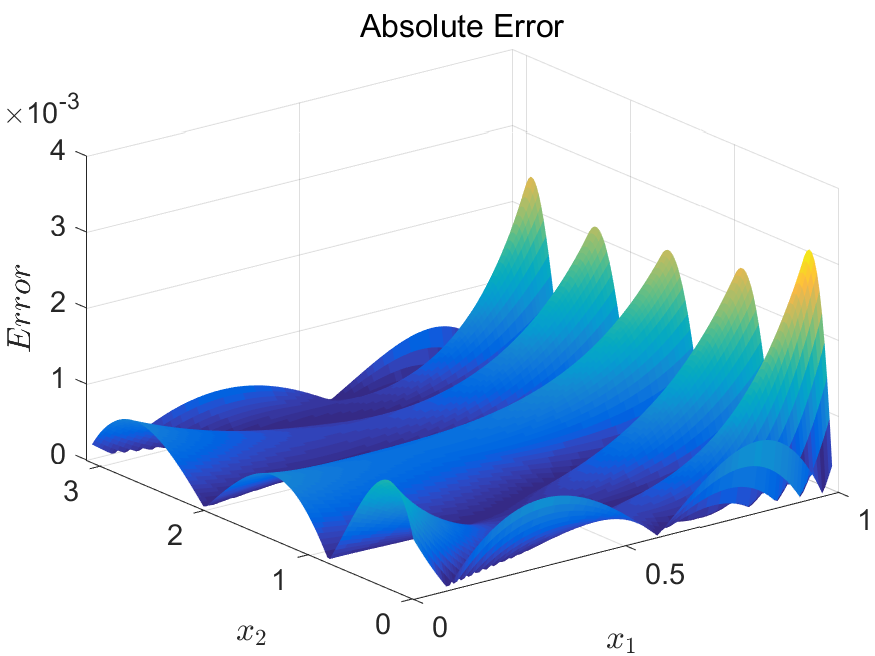}
	}
	\subfigure[Pointwise error for MIM]{
		\label{UABS2MIM_BiNavier2D_E1}
		\includegraphics[scale=0.34]{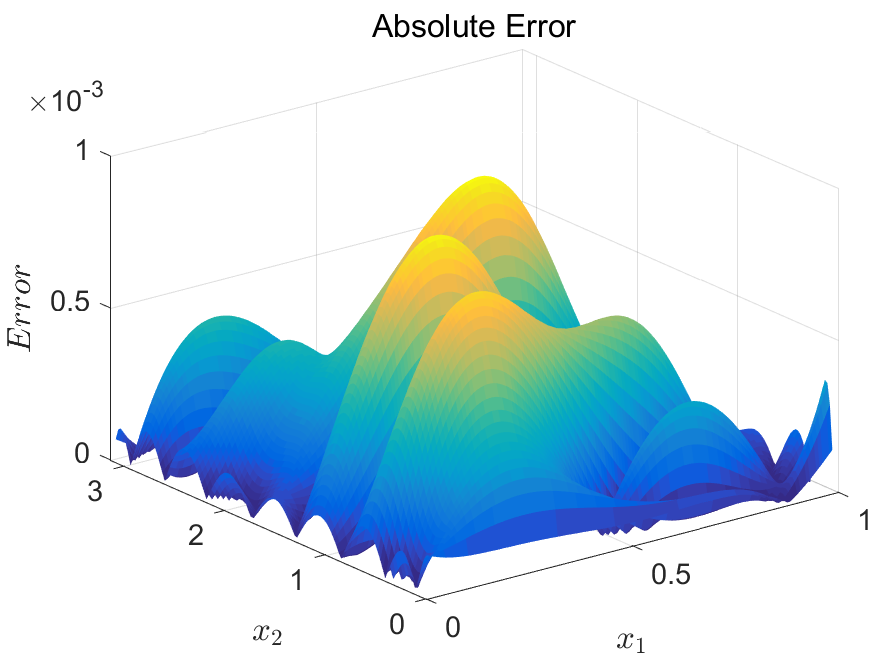}
	}
	\subfigure[Pointwise error for FCPINN]{
		\label{UABS2CPINN_BiNavier2D_E1}
		\includegraphics[scale=0.34]{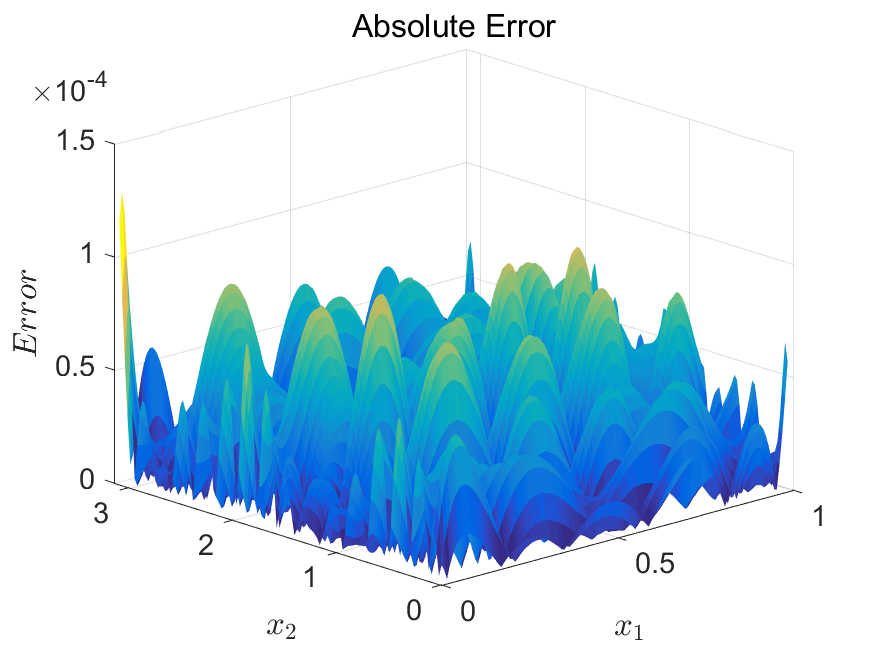}
	}
	\subfigure[REL]{
		\label{REL2BiNavier2D_E1}
		\includegraphics[scale=0.275]{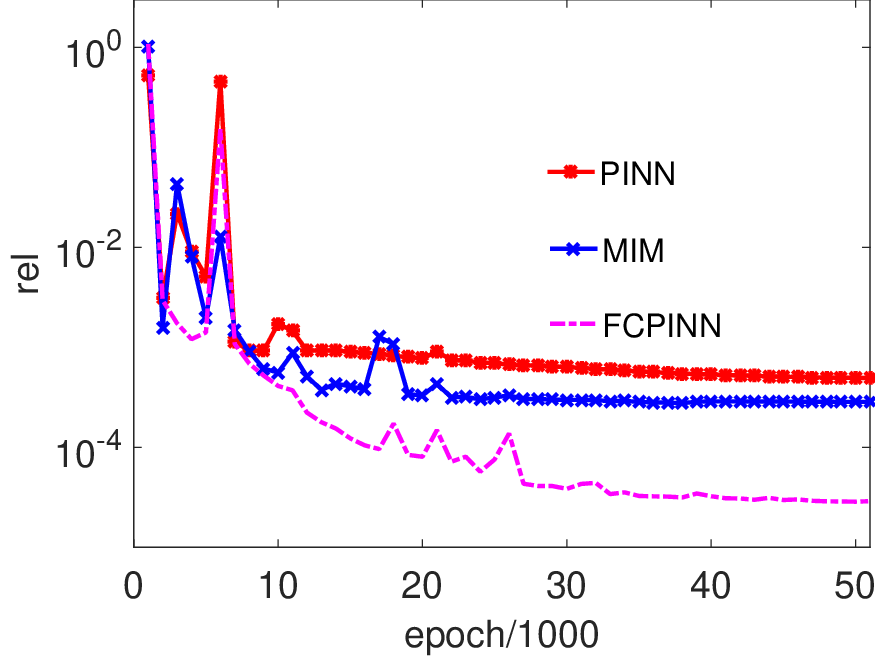}
	}
	\label{BiNavier_2DE1_U}
	\caption{The test results of solution for Example \ref{BiNavier2D_E01}.}
\end{figure}

\begin{figure}[H]
	\centering
	\subfigure[Exact solution]{
		\label{VExact2BiNavier2D_E1}
		\includegraphics[scale=0.325]{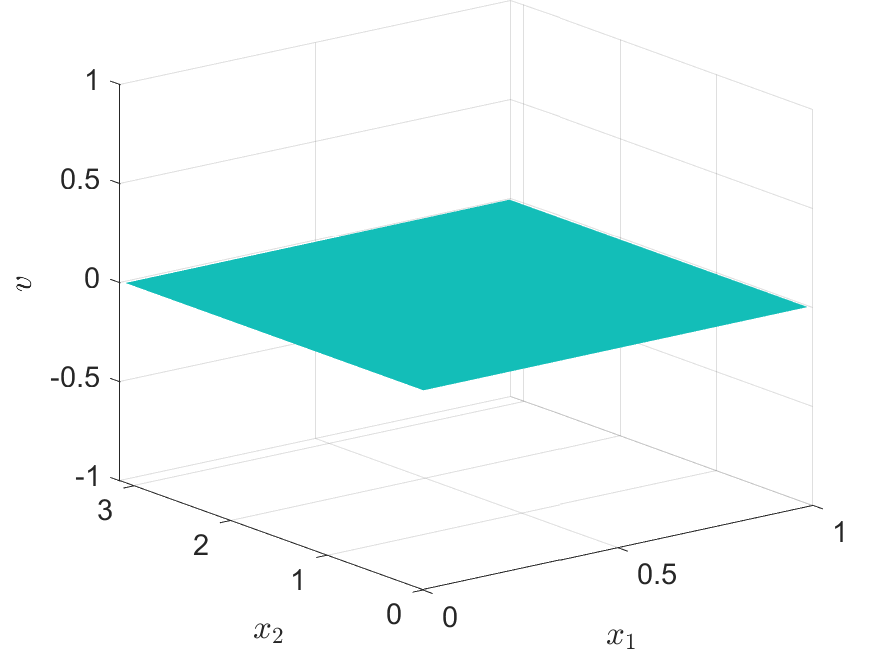}
	}
	\subfigure[Pointwise error for MIM]{
		\label{VABS2MIM_BiNavier2D_E1}
		\includegraphics[scale=0.325]{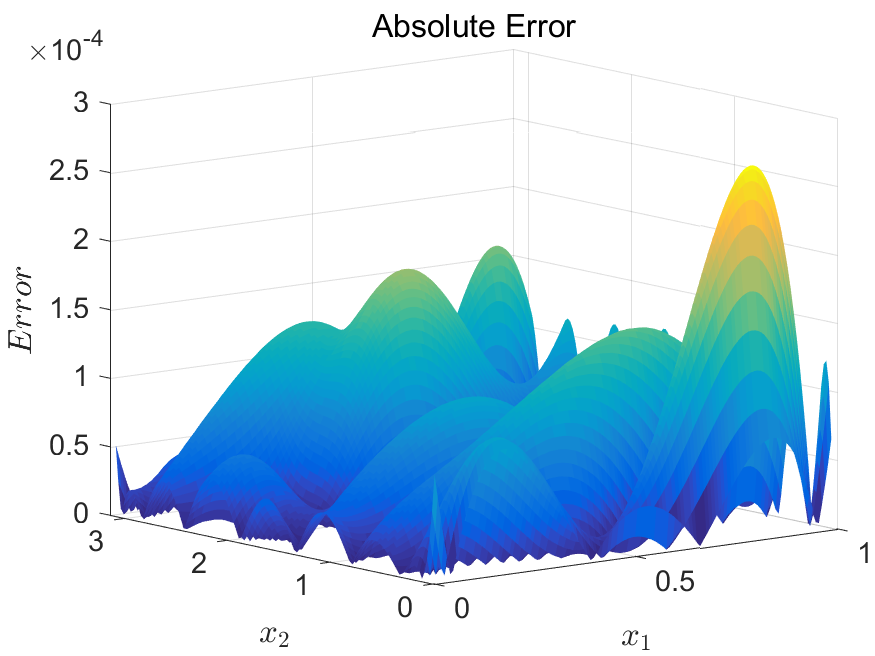}
	}
	\subfigure[Pointwise error for FCPINN]{
		\label{VABS2CPINN_BiNavier2D_E1}
		\includegraphics[scale=0.325]{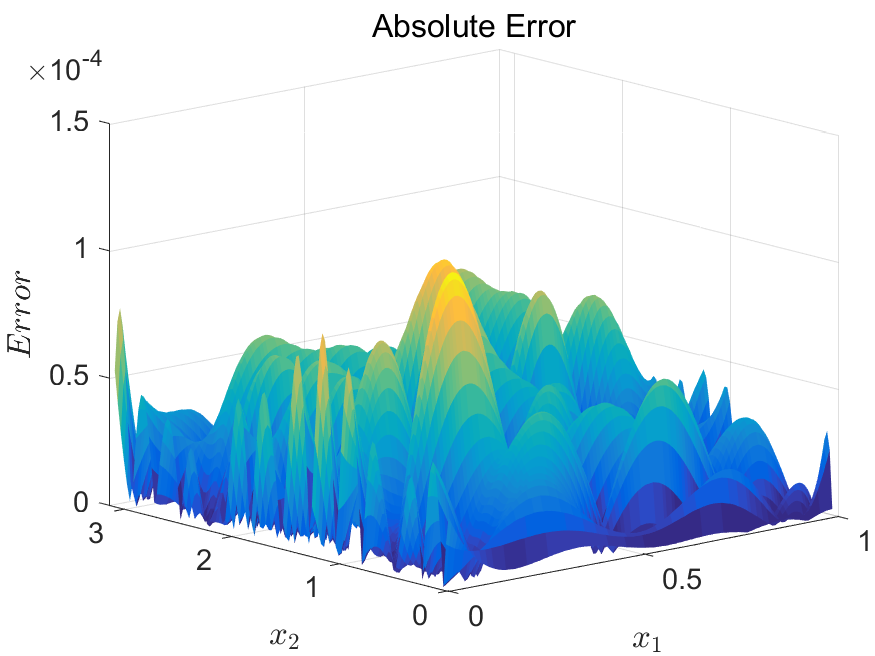}
	}
	\label{BiNavier_2DE1_V}
	\caption{The test results of $v$ for Example \ref{2D_Dirichlet_E1}.}
\end{figure}

\begin{example}\label{BiNavier2D_E2}
We consider a nonlinear case of the biharmonic equation~\eqref{eq:biharmonic} on a porous domain $\Omega$ derived from the two-dimensional domain $[-\pi,\pi]\times[-\pi,\pi]$, with the force term defined as $f(x_1,x_2;u, \Delta u) = \tilde{f}(x_1,x_2)-\Delta u -u$. The exact solution is given by
\begin{equation*}
u(x_1,x_2)=\sin(x_1)\sin(x_2),
\end{equation*}
such that $f(x_1,x_2) = 3\sin(x_1)\sin(x_2)$ and $k(x_1,x_2)=-2\sin(x_1)\sin(x_2)$. The boundary function $g(x_1, x_2)$ on $\partial \Omega$ can be derived directly from the exact solution. The network configuration for this problem is the same as that used in the linear case.
\end{example}

In this example, the network configuration remains the same as in the linear case, so we omit further details. For nonlinear problems with Navier boundary conditions, the performance of the three models follows a similar trend to that in linear cases: in terms of residual error, MIM performs worse than PINN, while FCPINN achieves the highest accuracy. Regarding computation time, MIM is the fastest, FCPINN is slightly slower, and PINN takes the longest time. Due to the weakened boundary constraints inherent in Navier boundary conditions, MIM loses many features of the original function in nonlinear problems, resulting in lower accuracy. In contrast, PINN and FCPINN, benefiting from the inherent complexity of DNN methods, are better suited for problems with implicit constraints. Additionally, in this case, the broader value range of the function highlights the effectiveness of Fourier activation layers in improving accuracy within a large-scale value space. This underscores FCPINN’s advantage in handling large-scale Navier boundary problems.

\begin{table}[!ht]
    \centering
    \caption{The REL and time consumed of the above three models on various domains for example \ref{BiNavier2D_E2} }
    \label{Table2BiNavier2D_E2}
    \begin{tabular}{|l|c|c|c|}
        \hline  
                  & PINN     & MIM                  &  FCPINN               \\ \hline
          REL     &$0.0011$  &$0.0067$ &$6.602\times10^{-4}$  \\ 
         Time(s)  &1656.033  &594.632              & 663.623               \\  \hline
    \end{tabular}
\end{table}

\begin{figure}[H]
    \centering
    \subfigure[Exact solution]{
        \label{UExact2BiNavier2D_E2}
        \includegraphics[scale=0.34]{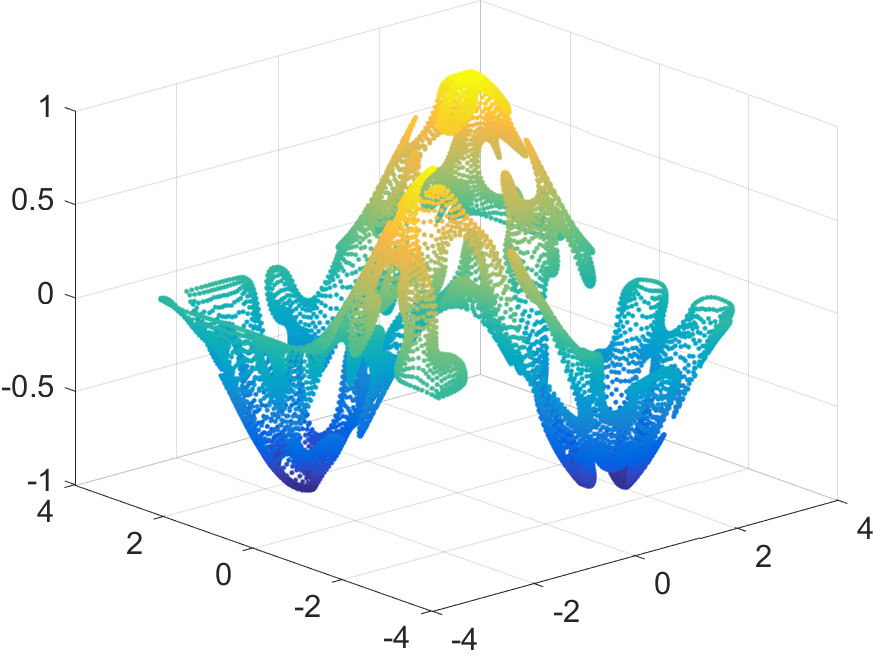}
    }
    \subfigure[Pointwise error for PINN]{
        \label{UABS2PINN_BiNavier2D_E2}
        \includegraphics[scale=0.34]{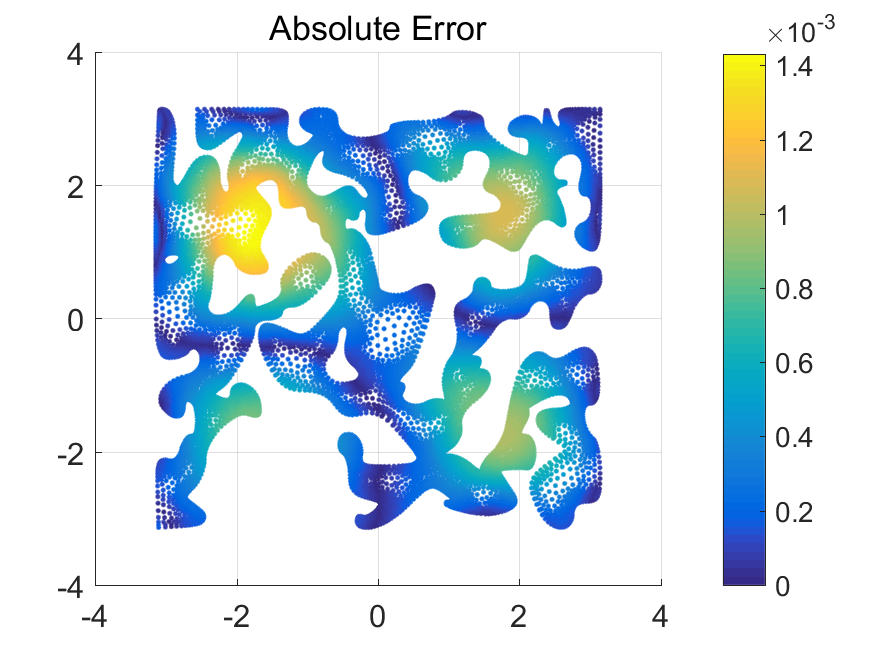}
    }
    \subfigure[Pointwise error for MIM]{
        \label{UABS2MIM_BiNavier2D_E2}
        \includegraphics[scale=0.34]{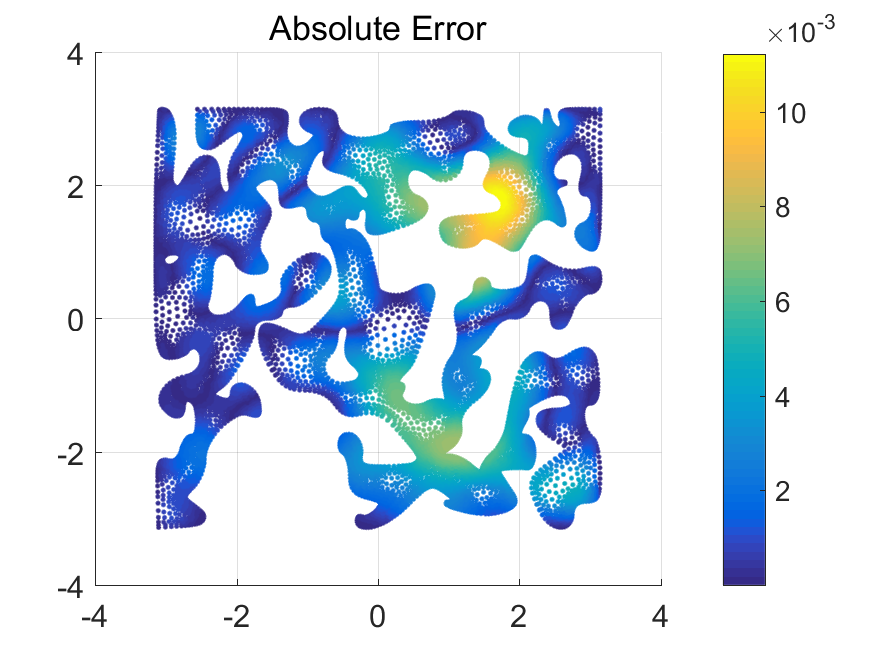}
    }
    \subfigure[Pointwise error for FCPINN]{
        \label{UABS2CPINN_BiNavier2D_E2}
        \includegraphics[scale=0.34]{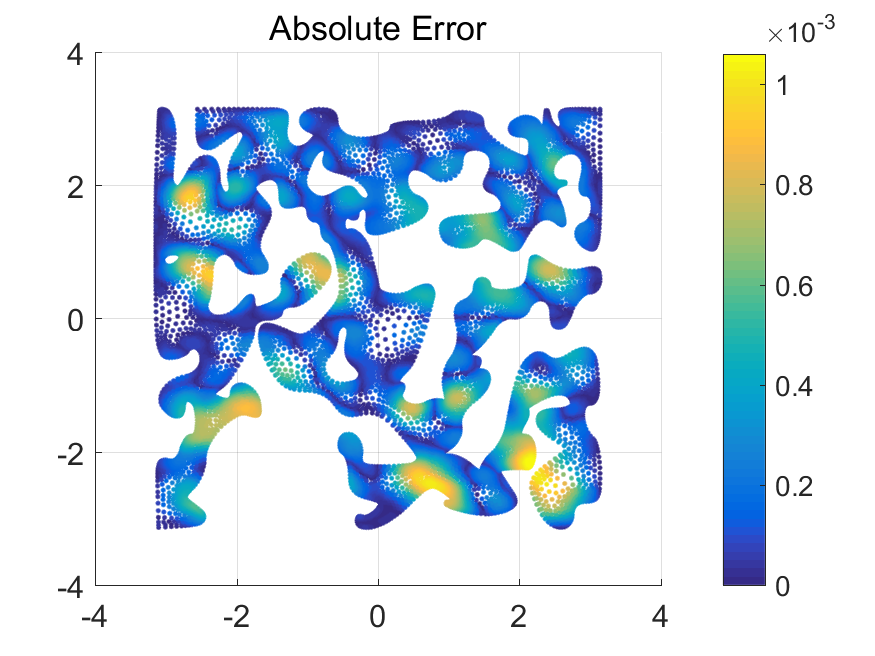}
    }
    \subfigure[REL]{
        \label{REL2BiNavier2D_E2}
        \includegraphics[scale=0.275]{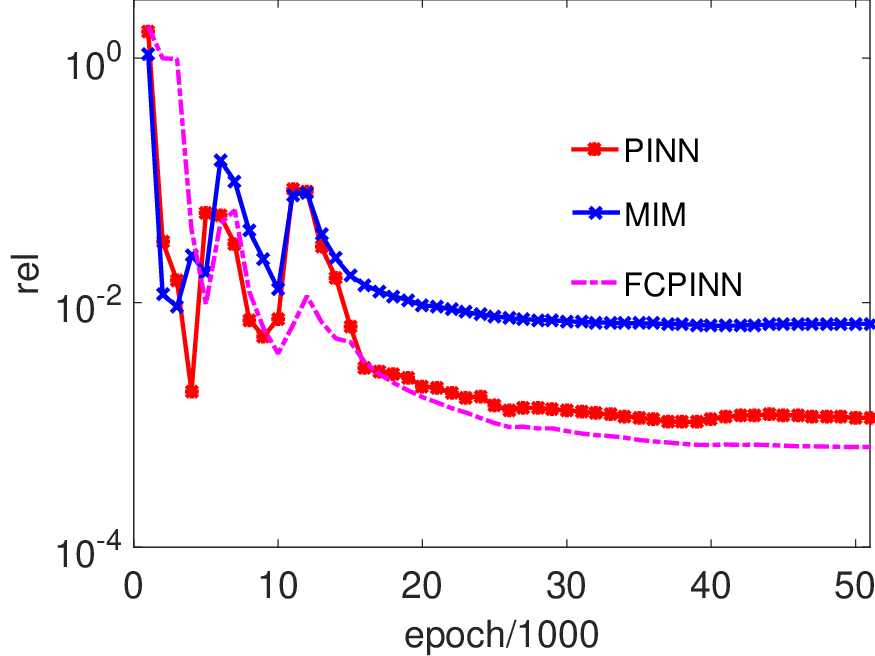}
    }
    \label{BiNavier_2DE2_U}
    \caption{The test results of solution for Example \ref{BiNavier2D_E01}.}
\end{figure}

\begin{figure}[H]
    \centering
    \subfigure[Exact solution]{
        \label{VExact2BiNavier2D_E2}
        \includegraphics[scale=0.325]{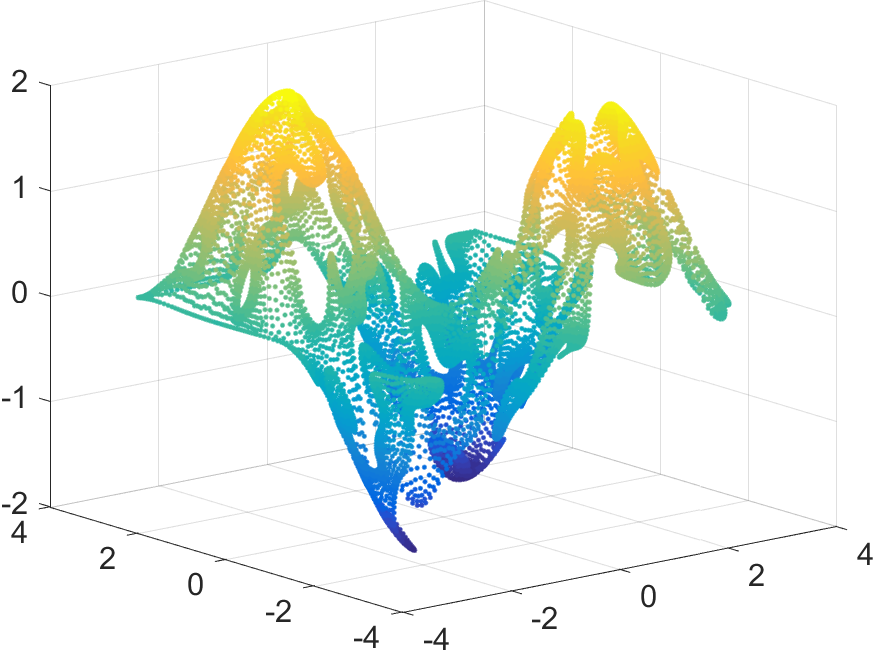}
    }
    \subfigure[Pointwise error for MIM]{
        \label{VABS2MIM_BiNavier2D_E2}
        \includegraphics[scale=0.325]{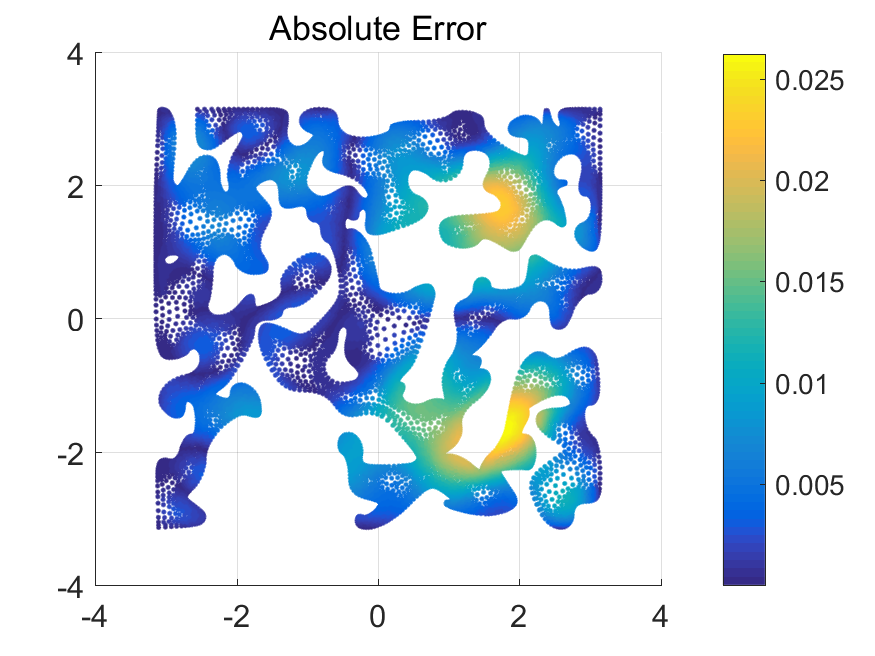}
    }
    \subfigure[Pointwise error for FCPINN]{
        \label{VABS2CPINN_BiNavier2D_E2}
        \includegraphics[scale=0.325]{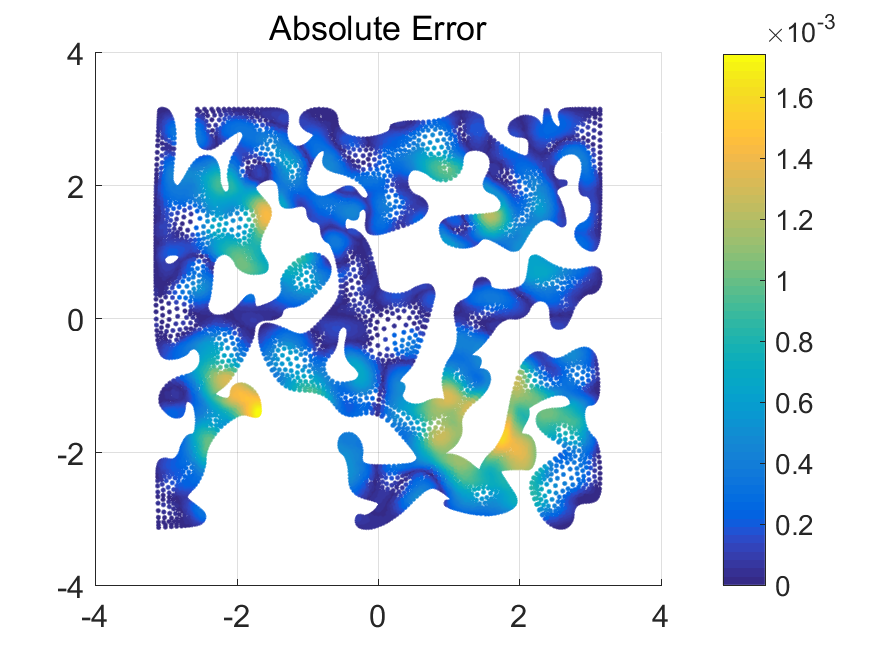}
    }
    \label{BiNavier_2DE2_V}
    \caption{The test results of $v$ for Example \ref{2D_Dirichlet_E1}.}
\end{figure}

\begin{example}\label{8D_Navier_E1}
In this example, we solve the biharmonic equation \eqref{eq:biharmonic} on an eight-dimensional unit-sized domain $\Omega=[-1,1]^8$. The exact solution is given by
\begin{equation*}
\begin{aligned}
    u(x_1,x_2,...,x_8)&=\sin(2\pi x_1)+\cos(2\pi x_2)+\sin(2\pi x_3)+\cos(2\pi x_4)\\
    &+\sin(2\pi x_5)+\cos(2\pi x_6)+\sin(2\pi x_7)+\cos(2\pi x_8),
\end{aligned}
\end{equation*}
so that the force term is
\begin{equation*}
f(x_1,x_2,...,x_8)=64\pi^4u(x_1,x_2,...,x_8).
\end{equation*}
The boundary conditions $g(x_1,x_2,...,x_8)$ and $h(x_1,x_2,...,x_8)$ on $\partial \Omega$ are easily derived as
\begin{equation*}
\begin{cases}
    \displaystyle g(x_1,x_2,...,x_8)=u(x_1,x_2,...,x_8),\\
    \displaystyle h(x_1,x_2,...,x_8)=-4\pi^2u(x_1,x_2,...,x_8).
\end{cases}
\end{equation*}
\end{example}

Due to the high computational resource requirements of PINN for 8D problems, we were unable to meet its resource demands. Therefore, only FCPINN and MIM were used to approximate the solution of the biharmonic equation in the 8D domain, with both models configured with network sizes of (40, 100, 80, 80). During training, we sampled 7000 points from the domain $\Omega$ and 500 points from the boundary $\partial \Omega$, while the test set included 1600 randomly distributed points within $\Omega$. The residual errors and convergence process are shown in the figure below.

As seen in Table 8, FCPINN maintains its strong performance for 8D problems, with a computation speed nearly matching that of MIM, but with significantly greater accuracy. However, compared to conventional 2D and 3D problems, the residual error is relatively larger, and the computation time is longer. These results highlight FCPINN's adaptability and stability in handling high-dimensional problems.

\begin{table}[!ht]
    \centering
    \caption{The REL and Time of the above MIM and FCPINN methods, for example, \ref{8D_Navier_E1} on various domains}
    \label{Table2BiNavier8D_E1}
    \begin{tabular}{|l|c|c|}
        \hline  
                  & MIM      &  FCPINN               \\ \hline
          REL     &0.0045    &$4.891\times 10^{-4}$  \\ 
         Time(s)  &3175.338               &3692.739               \\  \hline
    \end{tabular}
\end{table}

 \begin{figure}[H]
    \centering
    \includegraphics[scale=0.55]{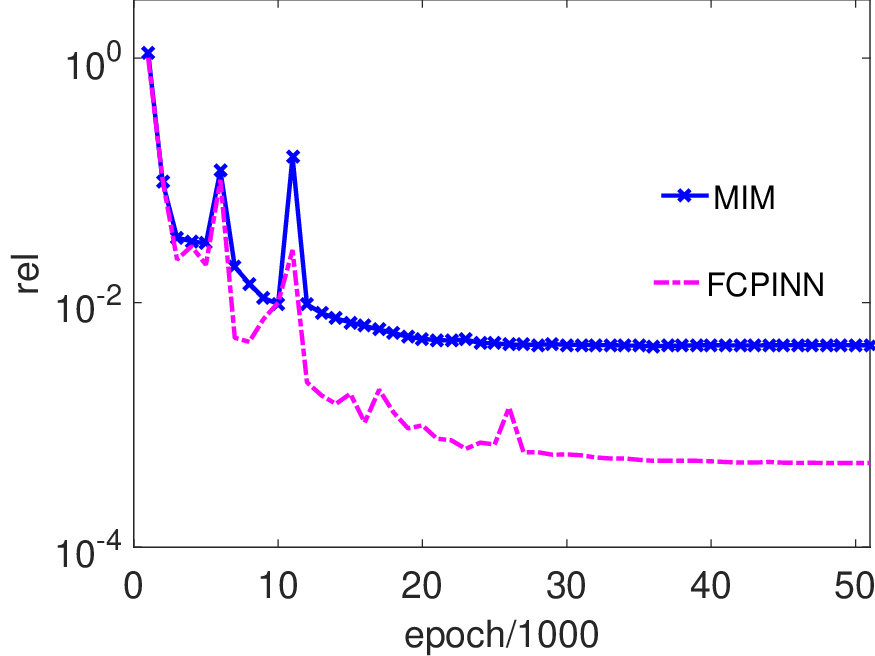}
    \label{fig:EightDim_Holes}
    \caption{REL for MIM and FCPINN methods to the biharmonic equation in the 8D domain}
 \end{figure}

\section{Conclusion}\label{sec:05}
In this paper, we have proposed a Coupled Physics-Informed Deep Neural Network (CPINN) framework to simplify the solution of complex biharmonic equations. Unlike traditional methods, the DNN approach—with its mesh-free nature and strong nonlinear approximation capabilities—can efficiently approximate solutions to biharmonic equations. Drawing inspiration from decoupling techniques in conventional methods, we introduce a set of auxiliary variables to facilitate computations, which proves effective for higher-order problems. To address the issue of spectral bias, we incorporate Fourier feature mapping as the activation function in the input layer, creating the FCPINN model. In the subsequent hidden layers, we use functions with strong regularization properties, such as sine. This blend of activation functions improves both the accuracy and convergence speed of the model. Numerical results confirm that this approach is feasible for solving \eqref{eq:biharmonic} in complex domains and high-dimensional spaces with high efficiency. In future work, we will further expand and refine the neural network model to tackle other higher-order PDEs.

\section*{Declaration of interests}
All authors confirm that they have no known financial conflicts or personal relationships that could have potentially influenced the work presented in this paper. 

\section*{Credit authorship contribution Statement}
Yujia Huang: Software, Data Curation, Visualization, Writing-Original Draft; Xi'an Li: Writing-Reviewing and Editing, Software, Visualization, Conceptualization, Formal analysis, Investigation; Jinran Wu: Supervision, Writing-Reviewing and Editing. All authors have reviewed and approved the final version of the manuscript for publication.

\section*{Acknowledgements}
This research was supported by the ``Chunhui'' Program Collaborative Scientific Research Project (202202004).


\end{document}